\documentclass[12pt]{article}
\usepackage[utf8]{inputenc}
\usepackage{amsmath}
\usepackage{amsfonts}
\usepackage{amsthm}
\usepackage{multicol}
\usepackage{enumerate}
\usepackage{varioref}
\usepackage[colorlinks=true, allcolors=black]{hyperref}
\usepackage{array}
\usepackage{cite}
\usepackage{ragged2e}
\usepackage{mathtools}
\usepackage{tabstackengine}
\usepackage{amssymb}
\usepackage{lipsum} 
\usepackage{graphicx} 
\usepackage{epstopdf}
\usepackage{color}

\usepackage{subcaption}
\usepackage{float}
\usepackage{textgreek}
\usepackage{babel}
\usepackage{tikz}
\usetikzlibrary{arrows.meta}
\usetikzlibrary{angles, quotes}
\usepackage[labelfont=bf]{caption}

\usepackage{microtype}

\newtheorem{theorem}{Theorem}
\newtheorem*{problem*}{Problem}

\newtheoremstyle{case}{}{}{}{}{}{:}{ }{}
\theoremstyle{case}
\newtheorem{case}{Case}
\theoremstyle{definition}
\newtheorem{definition}{Definition}

\usepackage{geometry}
\makeatletter
\let \@fnsymbol \@arabic
\makeatother
\title{Minimal surfaces over harmonic shears}
\author{
	Simran Bedi\thanks{Department of Mathematics, University of Delhi.\\
	sbedi@maths.du.ac.in}
	\and
	Sanjay Kumar \thanks{Deen Dyal Upadhyaya College, University of Delhi.\\
	skpant@ddu.du.ac.in}
}
\date{}
\begin{document}
	\maketitle	
	\begin{abstract}
		Harmonic mappings have long intrigued researchers due to their intrinsic connection with minimal surfaces. In this paper, 
		we investigate shearing of two distinct classes of univalent conformal mappings which are convex in horizontal direction with appropriate dilatations. Subsequently, we present a family of minimal surfaces constructed by lifting the harmonic mappings obtained through shear construction method given by Clunie-Sheil. Furthermore, we contribute to addressing an open problem partially, proposed by Boyd and Dorff, by identifying the resulting minimal surfaces for certain values of the parameters in one of the classes of mappings. Notably, this family of minimal surfaces transforms from the well-established Enneper surface to a Helicoid.
	\end{abstract}
	
	\section{Introduction}
The study of minimal surfaces has fascinated differential geometers, not solely because of the exotic structures that they exhibit but due to their inherent relationship with harmonic mappings. Results and properties of harmonic univalent functions are generally used to investigate the properties of minimal surfaces. The reason why we emphasize their univalency is because, geometrically, it means that the image of this map does not show any overlapping or self-intersection. Further, when we lift this univalent map from the complex plane to $\mathbb{R}^3$, the resulting surface inherently does not self-intersect itself, called a “minimal graph”. Construction of such rather special functions is an intricate task.
	\par One of the approach is to use the shearing  technique, introduced by Clunie and sheil-small in 1984 \cite{clunie1984harmonic}, whose fundamental building blocks are: a conformal mapping and a dilatation. It proved to be  a revolutionary paper in the study of harmonic mapping which gave rise to a lot of research in this field. Another revolutionary discovery was the Weierstrass-Enneper representation, which allows us to take a harmonic univalent function with an appropriate dilatation and lift it to a minimal graph. Several  researchers have used this technique \cite{jiang2019univalent}, however, it is often difficult to identify the resulting minimal graphs.
	 \par  One approach to recognize the minimal surface is to use a change of variables. In \cite{dorff2014family}, Dorff and Muir have constructed and identified a family of minimal graphs associated with the generalized Koebe function using this approach.
	 In \cite{boyd2014harmonic}, Boyd and Dorff proposed the following open problem, which is the prime object of our investigation:
	 \begin{problem*}
	 	Determine the minimal graphs formed by lifting harmonic univalent mappings in the paper “Gauss curvature estimates for minimal graphs” by Nowak and Woloszkiewicz \cite{nowak2011gauss}.
	 	\end{problem*}  We have contributed to answering this open problem partially, by recognizing the surfaces formed by lifting the harmonic maps for certain values, obtained through the shearing method applied on the conformal univalent mapping $F(z)=z/(1+c\,z+z^2)$, given in the aforementioned paper. We have also constructed a two-parameter family of harmonic mappings by shearing the same map with a different dilatation, $\omega(z)=z(z+a)/(1+az)$.
	 \par Also, We have chosen another family of conformal univalent maps which are convex in horizontal direction(CHD) given by,  $F(z)=z-(1/n^2)z^n ;n=2,3,4,\hdots$ \cite{dorff2012anamorphosis}. We have investigated its shearing for $\omega(z)=z^n$ and then have lifted the resulting harmonic maps to construct corresponding minimal graphs for even $n.$ 
	 We have plotted images of open unit disk under the conformal mappings, images of harmonic mappings, and the resulting minimal surfaces they give rise to.  The projections of these minimal surfaces can be observed in the images showing their corresponding harmonic mappings. 
	\section{Background}
	Now that we have some essential ideas, we need the background for the theory of Harmonic mappings. 
	This section will cover important definitions and theorems, which we will utilize throughout the rest of the paper.
	\par Let  $\mathbb{D}=\{z:|z|<1\}$ be the open unit disk in the complex plane. 
	\begin{definition}
		Let $S$ denotes the family of analytic functions on $\mathbb{D}$ that are normalized and univalent; that is,
	\begin{equation*} S=\{f:\mathbb{D} \rightarrow \mathbb{C} \mid f \mbox{ is analytic and univalent with } f(0)=0, f'(0)=1 \}.
		\end{equation*}
		\end{definition}
		The next theorem tells us that a harmonic function defined on $\mathbb{D}$  can be expressed in terms of analytic functions, known as its canonical decomposition \cite{dorff2012anamorphosis}.
	\begin{theorem}
		Let $f=u+\dot{\iota}v$ be a complex-valued harmonic function such that $f:D\rightarrow \mathbb{C}$, where $D$ is a simply-connected domain, then there exist analytic functions $h$ and $g$ such that $f=h+\bar {g}$.	\end{theorem}
		\begin{definition}
		Let $S_H^0$ be the family of complex-valued harmonic mappings on $\mathbb{D}$ that are univalent and normalized; that is, $$S_H^0=\{f:\mathbb{D}\rightarrow \mathbb{C} \mid f \mbox{ is harmonic, univalent with } h(0)=0, g(0)=0, h'(0)=1, g'(0)=0 \}.$$
	\end{definition}
		Note that the harmonic mapping $f=h+\bar{g}$  can also be expressed as $f=\operatorname{Re}{(h+g)}+\dot{\iota}\operatorname{Im}{(h-g)}$. By a result of Lewy \cite{lewy1936non}, $f=h+\bar{g}$ is locally univalent and sense-preserving in $\mathbb{D}$ if and only if the Jacobian of $f$, $J_F(z)=|h'(z)|^2-|g'(z)|^2>0$ for all $z\in\mathbb{D}.$ 
		The \textit{dilatation} of $f$ is defined as $\omega(z)=g'(z)/h'(z)$. Jacobian being positive is equivalent to the condition $|\omega(z)|<1$ for all $z\in \mathbb{D}$.
		\begin{definition} A domain $\Omega$ is said to be convex in the horizontal direction(CHD) if it has a connected intersection with every line parallel to the real axis. It is also called convex in the direction of real axis(CRA).
			\end{definition}
		The following theorem \cite{clunie1984harmonic} establishes the construction of a harmonic map with a specified dilatation which is a crucial tool for construction of minimal surfaces:
		\begin{theorem}[Clunie and Sheil-Small ]\label{2} A  harmonic map $f=h+\bar{g}$ locally univalent in $\mathbb{D}$ is a univalent mapping of $\mathbb{D}$ such that $f(\mathbb{D})$ is a CHD Domain if and only if , h-g is a conformal univalent mapping of $\mathbb{D}$ onto a CHD Domain.
			\end{theorem}
		This process is known as shearing or shear method and can be presented step by step as:
		\begin{enumerate}[(i)]
			\item Choose a  conformal univalent map $F$ and decompose it as $F=h-g.$
			\item Choose an appropriate dilatation $\omega$ and write it as $\omega=g'/h'.$
			\item Solve for $h$ and $g$ using equations obtained in the preceding two steps.
			\item Write $f=h+\bar{g}$, which is the required harmonic univalent map convex in horizontal direction obtained by \textit{shearing} a conformal map along parallel lines.
		\end{enumerate}
			Now, we lay some background about minimal surfaces before discussing their connection with harmonic mappings.
			\begin{definition}
				A  minimal surface is a surface whose mean curvature vanishes at all its points. 
				\end{definition}
			At each point, the bending upward in one direction is matched with the bending downward in the orthogonal direction. Every point on the surface is a saddle point. We can also define minimal surfaces in terms of first fundamental form, second fundamental form, area function and isothermal coordinates \cite{dorff2012soap}.
			 In this paper, we will come across two well-known minimal surfaces which are: Enneper's surface and helicoid.
		The helicoid is parametrized on $\mathbb{D} \setminus (-1,0)$ as:
		$$Y_{0}(z)=\left(\operatorname{Re}\left(z-\frac{1}{z}\right),\operatorname{Im}\left(z+\frac{1}{z}\right),2\,\operatorname{Im}(\log z)\right)$$
		Enneper's surface is parametrized on $\mathbb{D}$ as:
		$$Y_{2}(z)=\left(\operatorname{Re}\left(z-\frac{z^3}{3}\right),\operatorname{Im}\left(z+\frac{z^3}{3}\right),2\,\operatorname{Re}(-z^2)\right)$$
		\justifying Note that the negative sign in the third component function has the effect of reflecting the surface through xy-plane. Moreover, scalings by a factor, substitution by a $m\ddot{o}bius$ transformation and reflection across planes containing two axes does not affect the geometry of the surface.\\
	Karl Weierstrass and Alfred Enneper discovered a way to easily construct minimal surfaces by demonstrating their connection to harmonic mappings in the next theorem \cite{duren2004harmonic}. 
	\begin{theorem}
		[Weierstrass-Enneper representation]\label{3}
		Let $\Omega \subset \mathbb{C}$ be a simply connected domain containing the origin. If a minimal graph
		$$
		\{(u, v, F(u, v)): u+i v \in \Omega\}
		$$
		is parameterized by sense-preserving isothermal parameters $z=x+i y \in \mathbb{D}$, the projection onto its base plane defines $a$ harmonic mapping $w=u+i v=f(z)$ of $\mathbb{D}$ onto $\Omega$ whose dilatation is the square of an analytic function. Conversely, if $f=h+\bar{g}$ is a harmonic univalent mapping of $\mathbb{D}$ onto $\Omega$ with dilatation $\omega=g^{\prime} / h^{\prime}$ being the square of an analytic function, then with $z=x+i y \in \mathbb{D}$, the parameterization
		\begin{align*}
			\mathbf{X}(z)=\left(\operatorname{Re}\{h(z)+g(z)\}, \operatorname{Im}\{h(z)-g(z)\},2 \operatorname{Im}\left\{\int_0^z \sqrt{g^{\prime}(\zeta) h^{\prime}(\zeta)} d \zeta\right\}\right)
		\end{align*}
		defines a minimal graph whose projection into the complex plane is $f(\mathbb{D})$. Except for the choice of sign and an arbitrary additive constant in the third coordinate function, this is the only such surface.

	\end{theorem}

	\section{Shearing of one-slit and two-slit conformal mappings}
	\begin{figure}[H]
		\centering
		\begin{subfigure}[b]{0.35\textwidth}
			\centering
			\includegraphics[width=\textwidth]{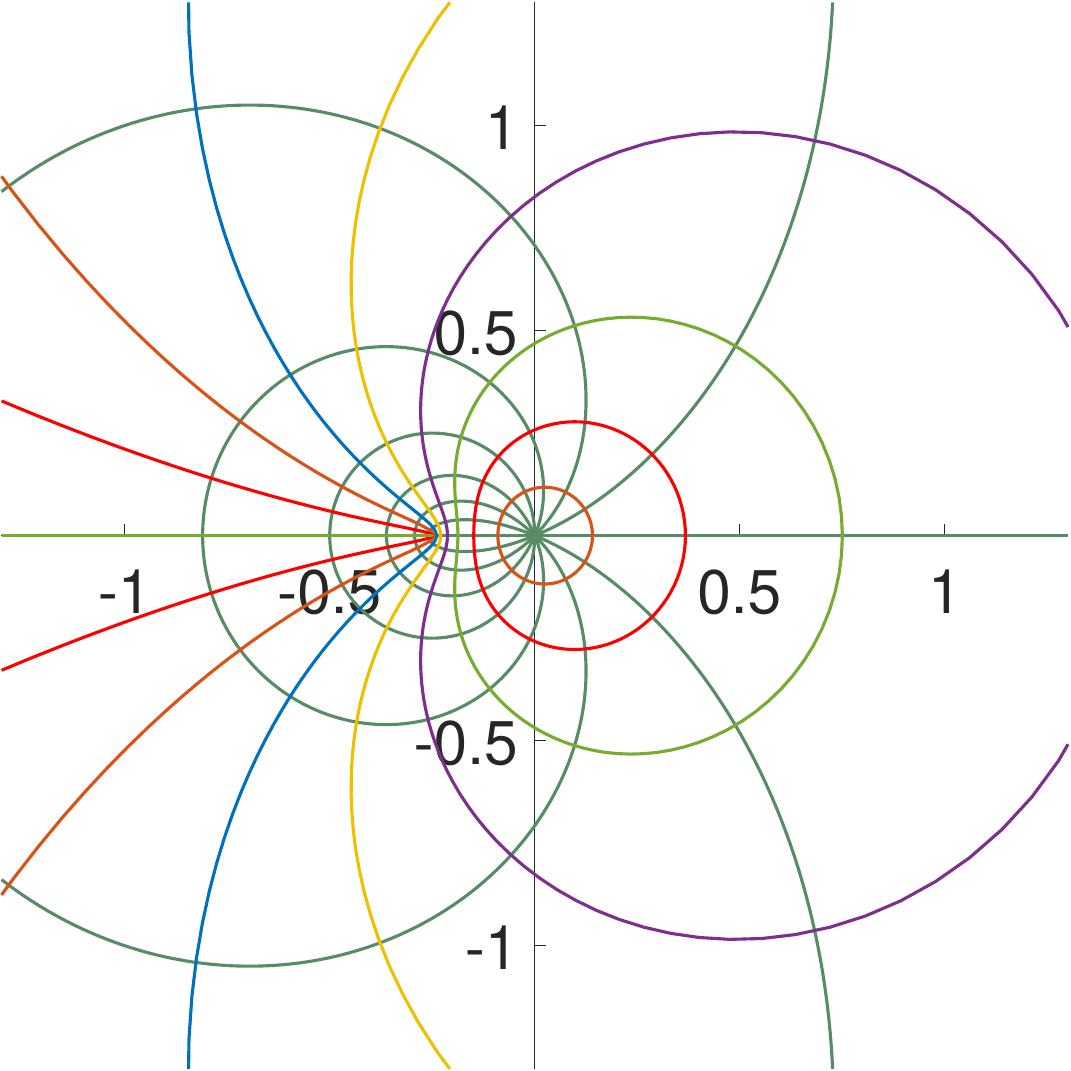}
			\caption{$c=-2$}
			\label{fig:c=-2}
		\end{subfigure}
		\begin{subfigure}[b]{0.35\textwidth}
			\centering
			\includegraphics[width=\textwidth]{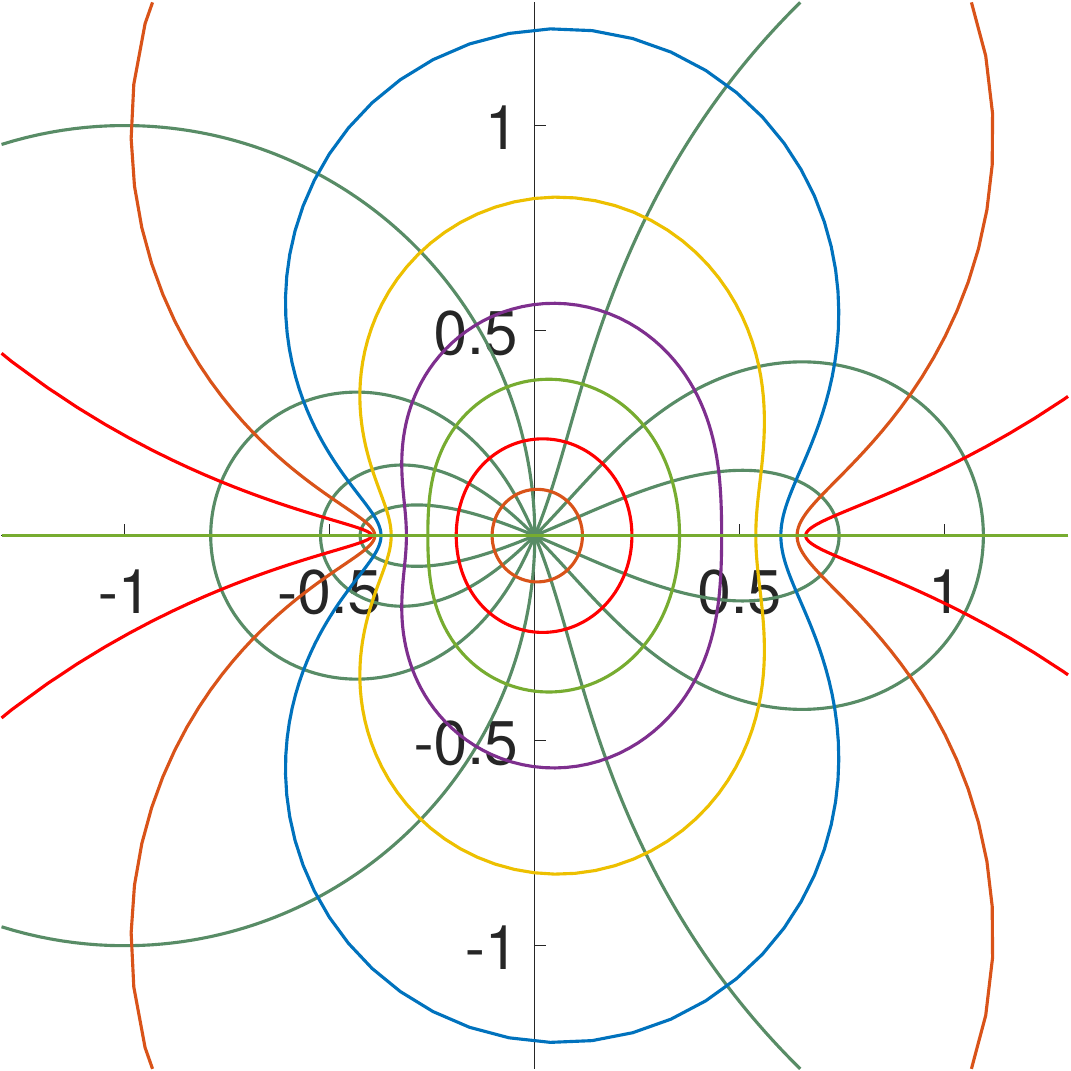}
			\caption{$c=0$}
			\label{fig:c=0}
		\end{subfigure}
		\begin{subfigure}[b]{0.35\textwidth}
			\centering
			\includegraphics[width=\textwidth]{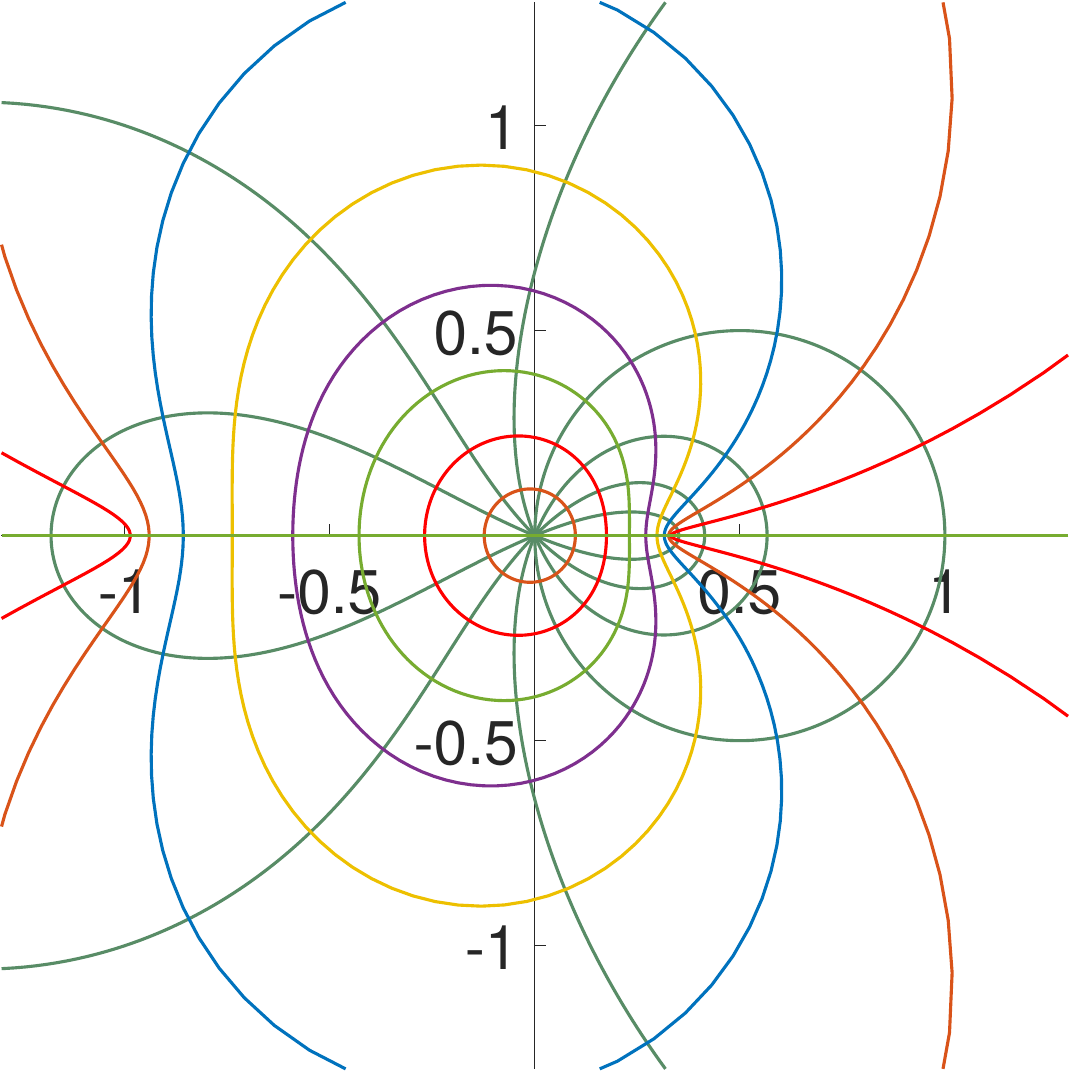}
			\caption{$c=1$}
			\label{fig:c=1}
		\end{subfigure}
		\begin{subfigure}[b]{0.35\textwidth}
			\centering
			\includegraphics[width=\textwidth]{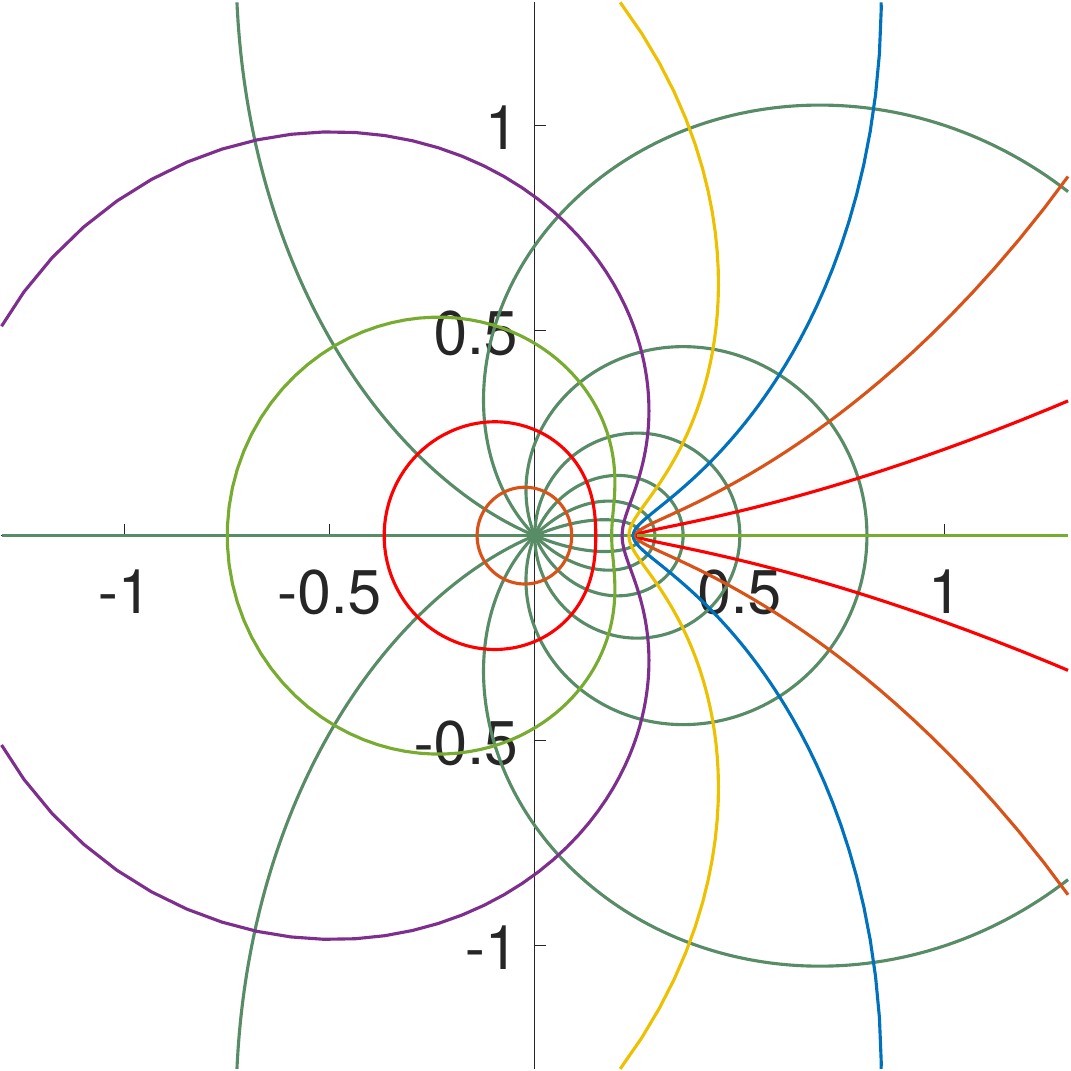}
			\caption{$c=2$}
			\label{fig:c=2}
		\end{subfigure}
		\caption{Image of  the conformal map $F_{c}(|z|=0.999)$}
		\label{fig:1}
	\end{figure}
		In this section, we apply shearing technique on a family of conformal mappings given in \cite{nowak2011gauss} denoted by $F_c(z)$ with two different dilatations $\omega_a(z)=z(z+a)(1+az)$ and $\omega(z)=z^n$. 
	We further lift the obtained two-parameter family of harmonic mappings $f_{c,a}$ from complex plane to minimal surfaces in $\mathbb{R}^3$ for $a=0$ and even $n$, using Theorem \ref{3}.
	Consider the family of univalent conformal map $F_c$ defined by,
	\begin{align}\label{eqn: 1}
		F_{c}(z)=\dfrac{z}{1+cz+z^2}
	\end{align} 
	where $c\in(-2,2)$.
	It maps the unit disk univalently onto a domain convex in the horizontal direction \cite{livingston1997univalent}.\\
	For special case of $c=-2$, the image of  under $F_{c}(\mathbb{D})$ is a single slit domain like Koebe domain depicted in Fig: \ref{fig:c=-2}, which is entire complex plane except for a slit along the negative real axis represented as: $$\mathbb{C} \setminus \{x: x\in(-\infty,-1/4)\}$$ and similarly a slit along positive real axis for $c=2$ shown in Fig: \ref{fig:c=2} given by: $$\mathbb{C} \setminus \{x: x\in(1/4,\infty)\}.$$ 
	For $c=0$, it maps the unit disk onto a two slit domain represented in Fig: \ref{fig:c=0}, that is, entire complex  plane except for two half-lines given by:
	$$\mathbb{C} \setminus \{x: x\in(-\infty,-1/2)\cup (1/2,\infty)\}.$$
	Let us consider the harmonic shear of $F_{c}$ defined in (\ref{eqn: 1}), with the dilatation
	$$\omega_{a}(z)=\frac{z(z+a)}{(1+az)},$$ such that $-1\leq a\leq 1$. In particular,  for $a=-1$, $a=1$ and $a=0$, $\,\omega(z)$ takes value $-z, z$ and $z^2$ respectively. This requires solution of the pair of differential equations given by:
	\begin{align*}
		h_{c,a}'(z)-g_{c,a}'(z)=F_c'(z) \,\,\,\mbox{and}\,\,\, g_{c,a}'(z)=z\frac{(z+a)}{(1+az)} h_{c,a}'(z).
	\end{align*}
	After straightforward but tedious computations for normalized $h_{c,a}$ and $g_{c,a}$, the  solution is expressed as:
	\begin{align}\label{eqn: 2}
		h_{c,a}(z)=-\frac{p(z)}{q(z)},
	\end{align}
	\begin{align*}
			\text{where} \,\,\,\,p(z)=&16 \sigma_2+16 \sigma_1-4c^2 \sigma_2+16z^2 \sigma_2-4c^2\sigma_1+16z^2\sigma_1+2z\sigma_3\sigma_4-4c^2z^2\sigma_2\\
		&+2ac^3\sigma_1-4c^3z\sigma_1-8ac\sigma_2+16c z\sigma_2-4c^2z^2\sigma_1+2ac^3 \sigma_2-4c^3 z\sigma_2\\
		& -8ac\sigma_1 +16cz\sigma_1 +2ac^3 z^2 \sigma_1 -8acz^2 \sigma_2 -8ac^2 z\sigma_2\\
		&+2ac^4 z\sigma_2 +2ac^3 z^2\sigma_2 -8acz^2\sigma_1 -8ac^2 z\sigma_1 \\
		&+2ac^4z\sigma_1 +2az^2 \sigma_3 \sigma_4 -cz^2 \sigma_3 \sigma_4\\
		&-c^2z\sigma_3  \sigma_4+acz\sigma_3 \sigma_4 ,\\ 
		\\  
		q(z)=&(c^2-4)\sigma_3\, \sigma_4 (1+c z+z^2),
		\end{align*}
		such that$\,\,
		\sigma_1 =\mathrm{atanh}\left(\frac{{\left(c^2 -4\right)}\,{\left(c+2\,z\right)}}{\sigma_3 \,\sigma_4 }\right),
		\;\sigma_2 =\mathrm{atanh}\left(\frac{4\,c-c^3 }{\sigma_3 \,\sigma_4 }\right),
		\;\sigma_3 ={{\left(c-2\right)}}^{3/2} \text{and}
		\;\sigma_4 ={{\left(c+2\right)}}^{3/2}.$\\
	Similary, solving for $g_{c,a}(z)$ gives, \begin{equation}\label{eqn: 3}
		g_{c,a}(z)=-r(z)/q(z),
	\end{equation}
	\begin{align*}
			\text{where}\,\,\,\,\,\, r(z)&=16\sigma_2 +16\sigma_1 -4c^2 \sigma_2 +16z^2 \sigma_2 -4c^2 \sigma_1 +16z^2 \sigma_1 -2z\sigma_3 \sigma_4 -4c^2 z^2 \sigma_2 \\
		&+2ac^3 \sigma_1 -4c^3 z\sigma_1 -8ac\sigma_2 +16cz\sigma_2 -4c^2 z^2 \sigma_1 +2ac^3 \sigma_2 -4c^3 z\sigma_2 \\
		&-8ac\sigma_1+16cz\sigma_1 +2ac^3 z^2 \sigma_1 -8acz^2 \sigma_2 -8ac^2 z\sigma_2 \\
		&+2ac^4 z\sigma_2 +2ac^3 z^2 \sigma_2 -8acz^2 \sigma_1-8ac^2 z\sigma_1 \\
		&+2ac^4 z\sigma_1 +2az^2 \sigma_3 \sigma_4 -cz^2 \sigma_3 \sigma_4 \\
		&+a\,c\,z\sigma_3 \sigma_4 
	\end{align*}
	and $q(z),\, \sigma_1,\,\sigma_2,\,\sigma_3$ are as above.
	\par So, $f_{c,a}(z)=h_{c,a}(z)+\overline{g_{c,a}(z)}$ is the desired map corresponding to conformal map $F_{c,a}$. By Theorem \ref{2}, we have $f_{c,a} \in S_{H}^0$ and is convex in direction of real axis. So, $f_{c,a}$ is the required harmonic map  obtained by applying shear construction.
   \\Next we discuss special cases by taking $c=-2$ and $c=2$ in (\ref{eqn: 1}). By shearing
   $$F_{-2}(z)=\frac{z}{(1-z)^2}\,\,\text{with dilatation} \,\,\, \omega_a(z)=\frac{z(z+a)}{1+az},\,\,\,\,$$
   \begin{figure}[H]
   	\centering
   	\begin{subfigure}[b]{0.32\textwidth}
   		\centering
   		\includegraphics[width=\textwidth]{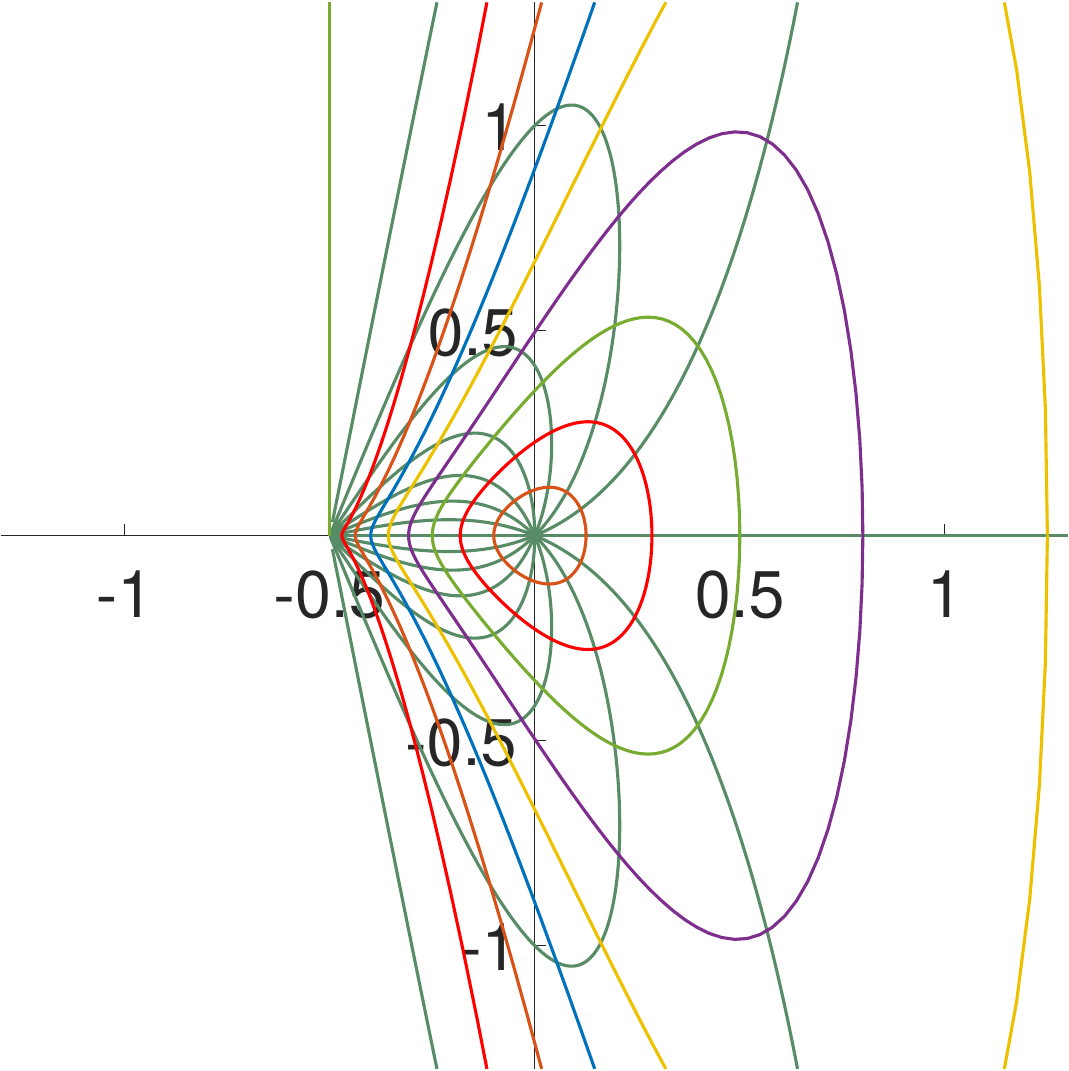}
   		\caption{$c=-2, a=-1$}
   		\label{fig: 2a}
   	\end{subfigure}
   	\begin{subfigure}[b]{0.32\textwidth}
   		\centering
   		\includegraphics[width=\textwidth]{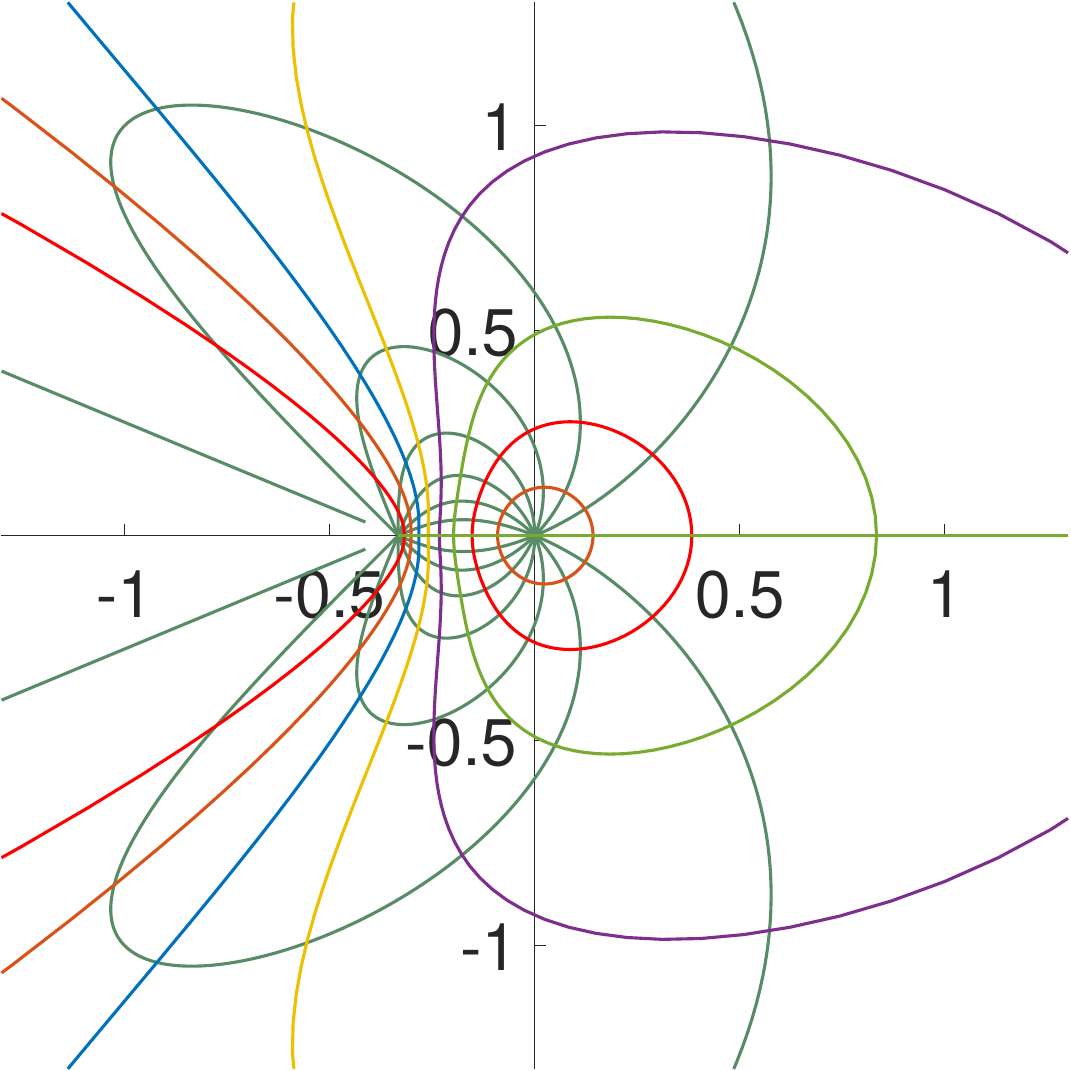}
   		\caption{$c=-2, a=0$}
   		\label{fig: 2b}
   	\end{subfigure}
   	\begin{subfigure}[b]{0.32\textwidth}
   		\centering
   		\includegraphics[width=\textwidth]{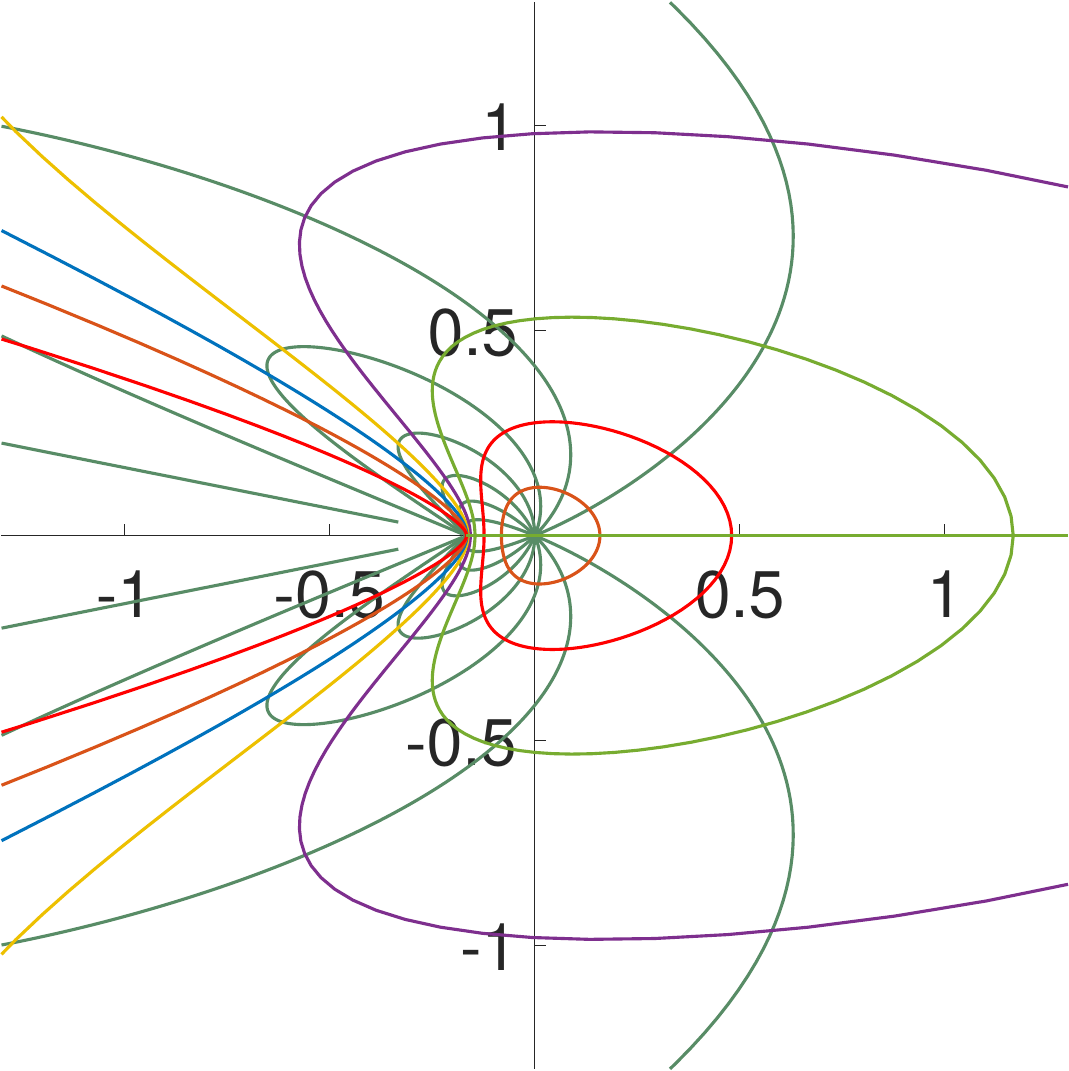}
   		\caption{$c=-2, a=1$}
   		\label{fig: 2c}
   	\end{subfigure}
   	\begin{subfigure}[b]{0.32\textwidth}
   		\centering
   		\includegraphics[width=\textwidth]{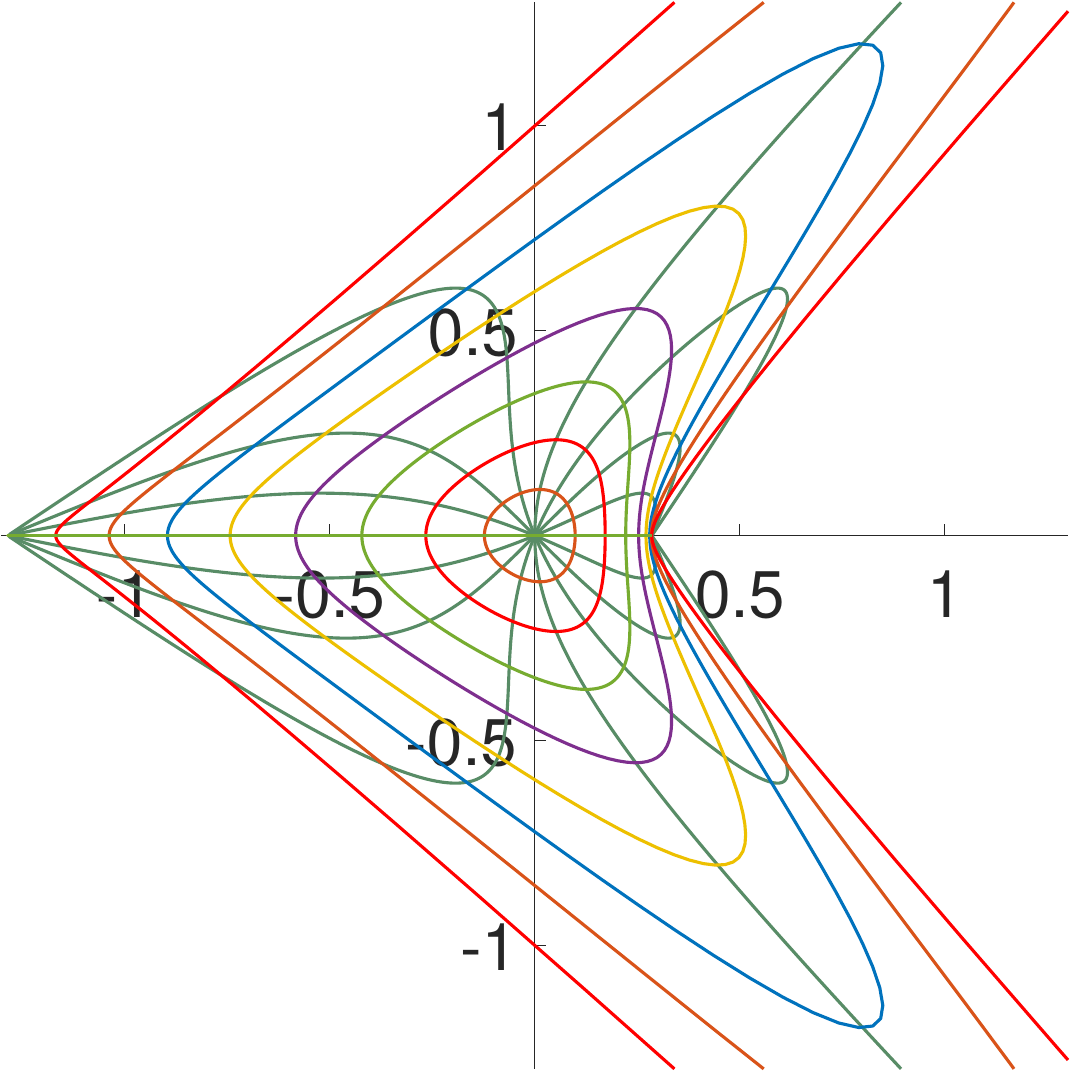}
   		\caption{$c=0, a=-1$}
   		\label{fig: 2d}
   	\end{subfigure}
   	\begin{subfigure}[b]{0.32\textwidth}
   		\centering
   		\includegraphics[width=\textwidth]{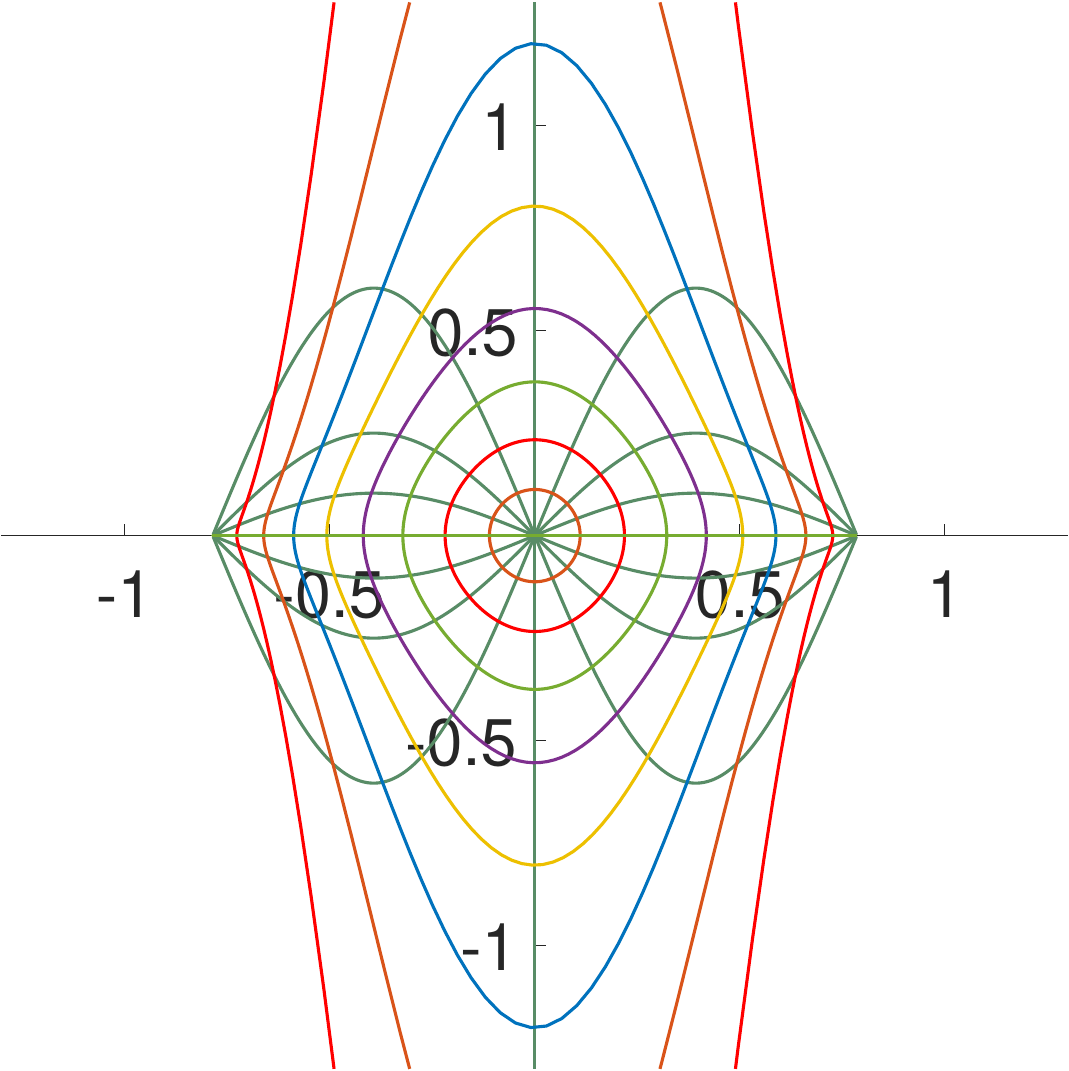}
   		\caption{$c=0, a=0$}
   		\label{fig: 2e}
   	\end{subfigure}
   	\begin{subfigure}[b]{0.32\textwidth}
   		\centering
   		\includegraphics[width=\textwidth]{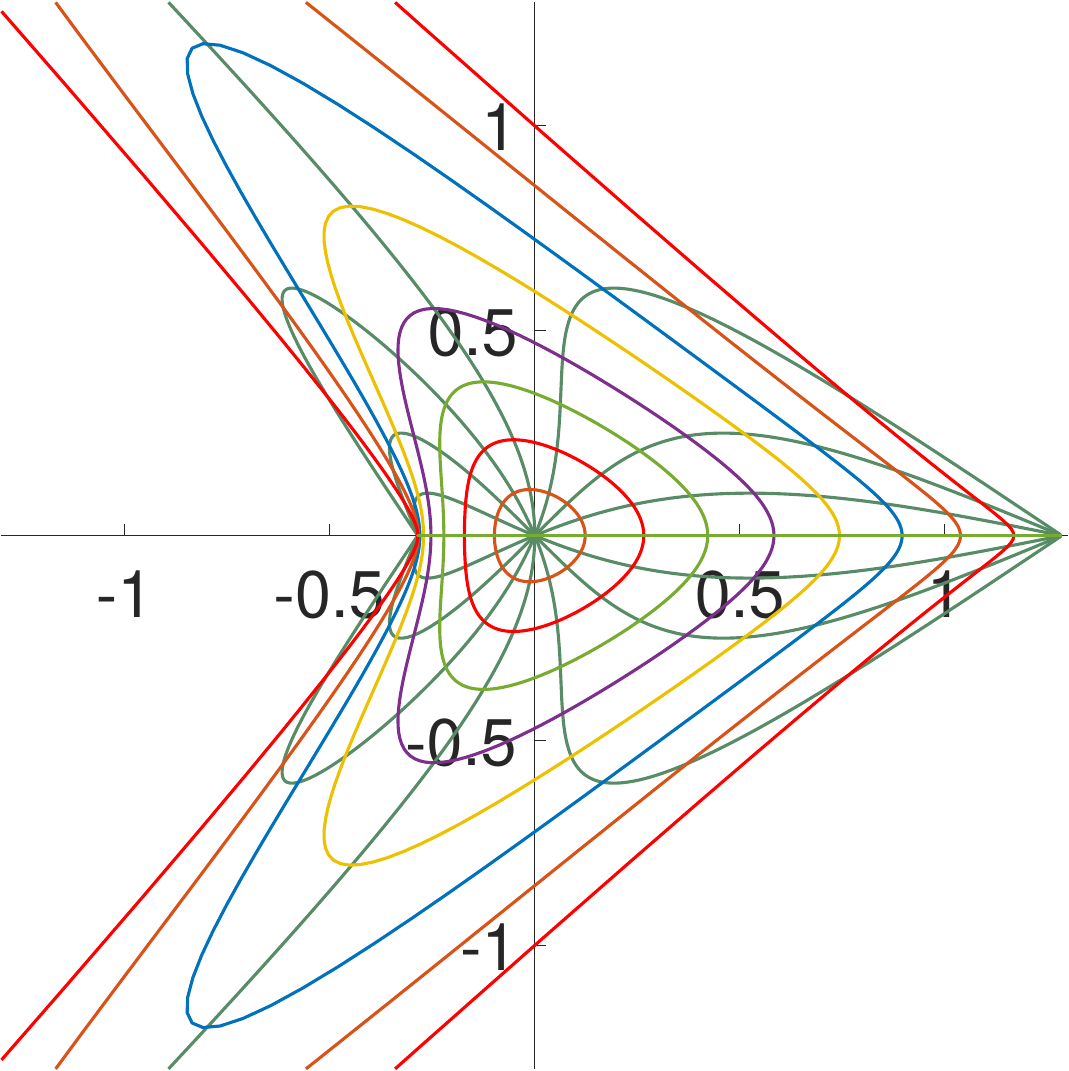}
   		\caption{$c=0, a=1$}
   		\label{fig: 2f}
   	\end{subfigure}
   	\begin{subfigure}[b]{0.32\textwidth}
   		\centering
   		\includegraphics[width=\textwidth]{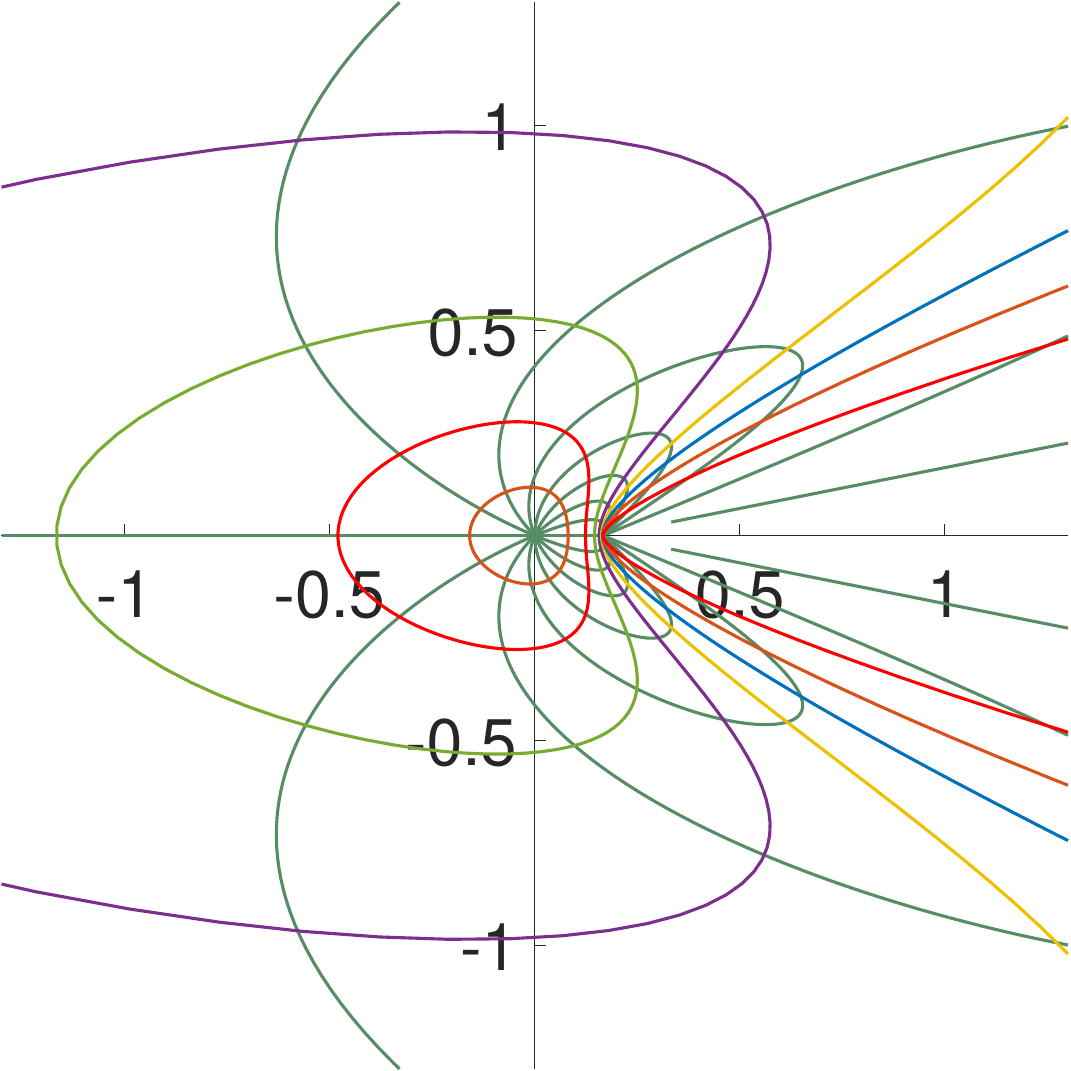}
   		\caption{$c=2, a=-1$}
   		\label{fig: 2g}
   	\end{subfigure}
   	\begin{subfigure}[b]{0.32\textwidth}
   		\centering
   		\includegraphics[width=\textwidth]{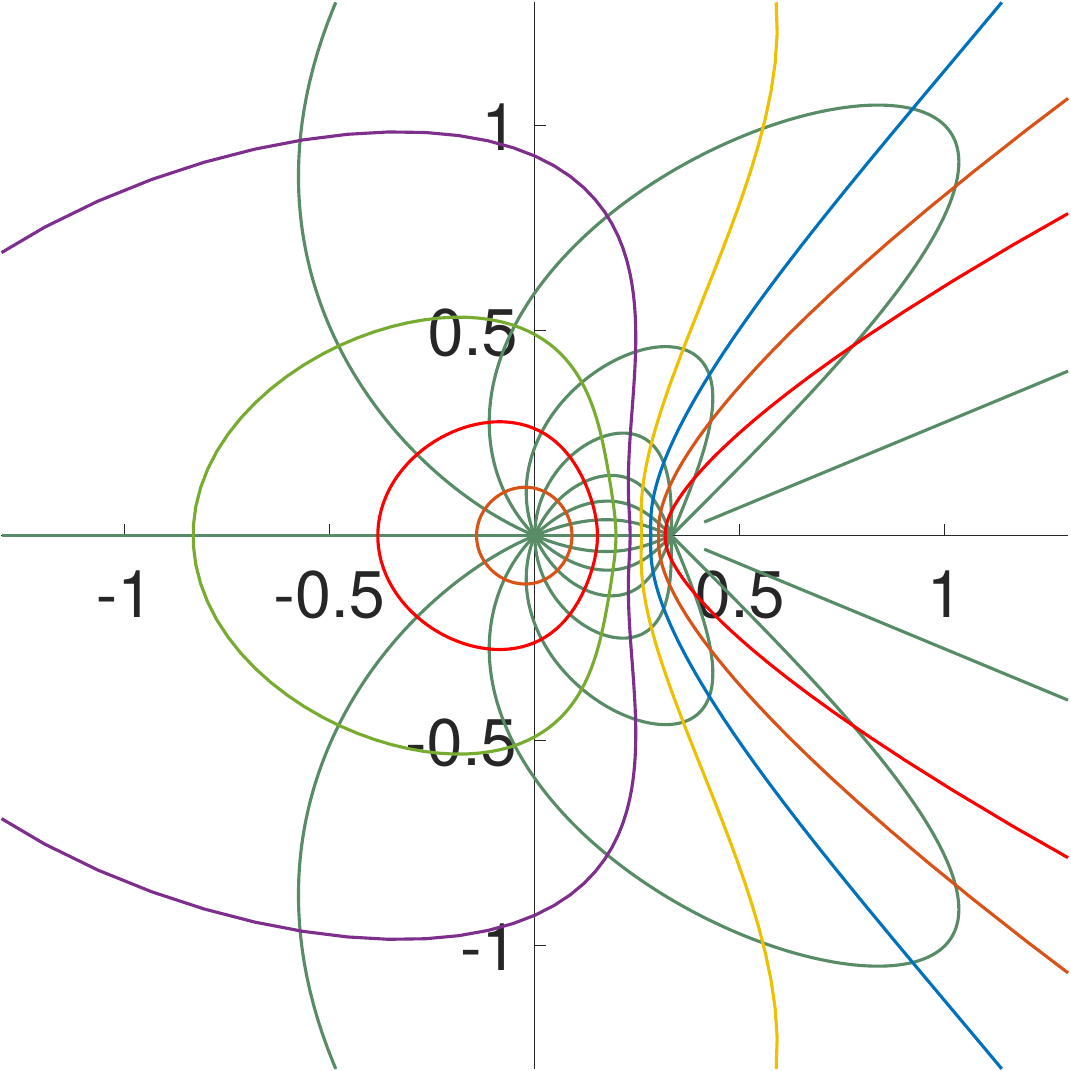}
   		\caption{$c=2, a=0$}
   		\label{fig: 2h}
   	\end{subfigure}
   	\begin{subfigure}[b]{0.32\textwidth}
   		\centering
   		\includegraphics[width=\textwidth]{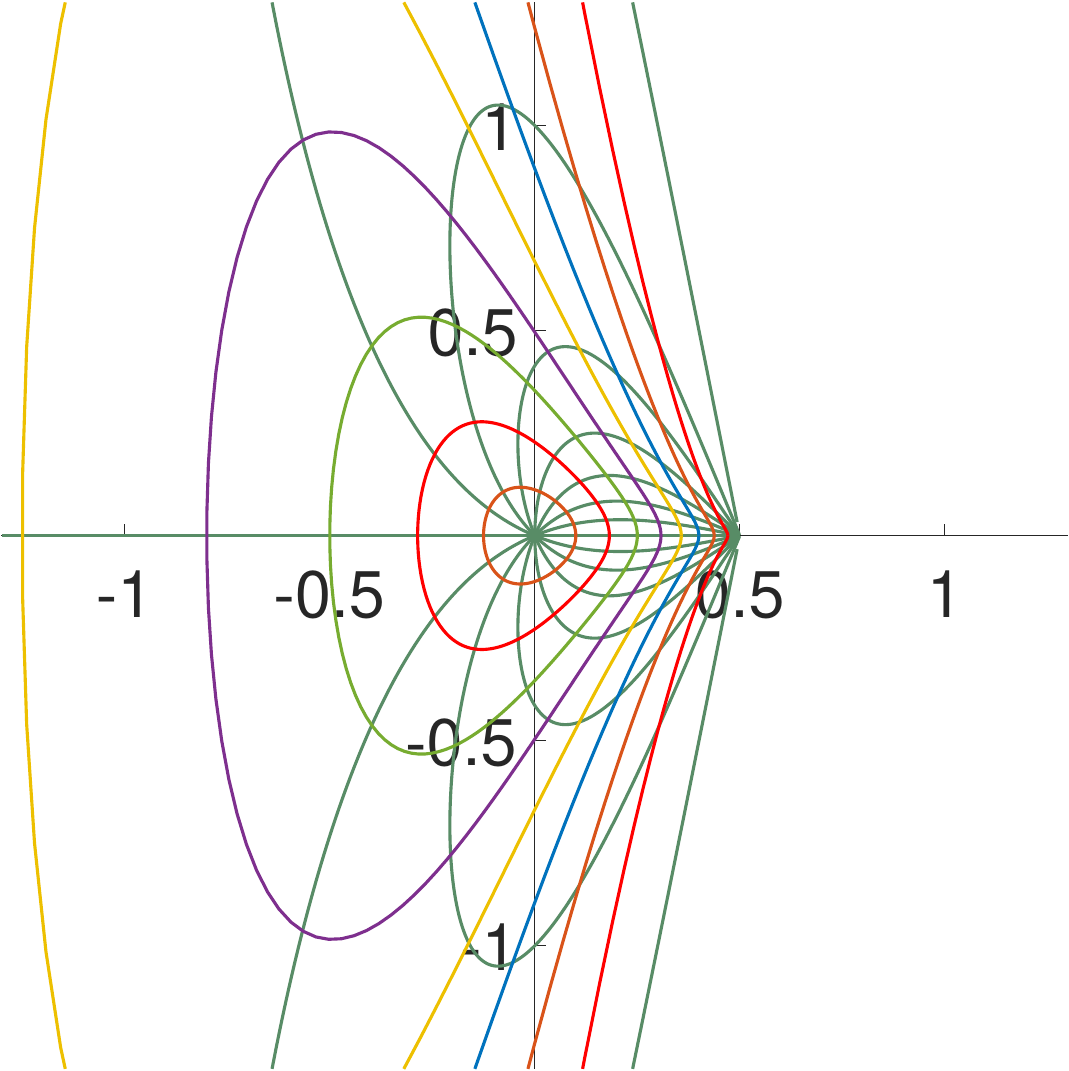}
   		\caption{$c=2, a=1$}
   		\label{fig: 2i}
   	\end{subfigure}
   	\caption{Images of harmonic shear $f_{c,a}(|z|=0.999)$ having dilatation $\omega(z)=\frac{z(z+a)}{1+az}$.}
   	\label{fig: 2}
   \end{figure}
   and solving for $h$ and $g$, we get $$h(z)=-\frac{6\,z+3\,a\,z^2 -a\,z^3 -6\,z^2 +2\,z^3 }{6\,{{\left(z-1\right)}}^3 } \,\,\,\text{and}\,\,\, g(z)=-\frac{3\,a\,z^2 -a\,z^3 +2\,z^3 }{6\,{{\left(z-1\right)}}^3 }$$
		So, the required map $f_{-2,a}$ can be expressed as $f_{-2,a}(z)=\operatorname{Re}\{h(z)+g(z)\}+\dot{\iota}\operatorname{Im}\{h(z)-g(z)\}$, i.e,
		\begin{align}\label{eqn: 4}
			f_{-2,a}(z)=\operatorname{Re}\left\{\frac{a}{3}-\frac{z^2+(a-1)z-\frac{a}{3}+\frac{2}{3}}{(z-1)^3}-\frac{2}{3}\right\}+\frac{\dot{\iota}}{4}\operatorname{Im}\left\{\left(\frac{1+z}{1-z}\right)^2-1\right\}
		\end{align}	
		By Theorem \ref{2}, we see that   $f_{-2,a}\in S_H^0.$ For $a=-1$, $f(\mathbb{D})$ is a right half plane. In order to study mapping properties of $f_{-2,a}$, let $w=(1+z)/(1-z)$, i.e, $z=(w-1)/(w+1)$, which leads us to:
		\begin{align*}
			f_{-2,a}(z)&=\operatorname{Re}\left\{\frac{1}{4}\left(\frac{1+a}{3}w^3+(1-a)w-\frac{2(2-a)}{3}\right)\right\}+\dot{\iota}\operatorname{Im}\left\{\frac{1}{4}(w^2-1)\right\}\\
			&=\frac{1}{4}\left(\frac{1+a}{3}(x^3-3xy^2)+(1-a)x-\frac{2(2-a)}{3}\right)+\dot{\iota}\frac{1}{2}xy, \,\, x>0.
		\end{align*}
		Note that each point $z\neq1$ on the unit circle is carried onto a point $w$ on the imaginary axis so that $x=0$ and $F_{-2,a}=-(2-a)/6.$ Similar discussion as in the case of harmonic koebe function K(z) in \cite{duren2004harmonic} shows that for $-1<a \leq 1$, $f(\mathbb{D})$ is a slit domain i.e, it maps onto the entire plane minus the interval $(-\infty,-(2-a)/6)$ on the negative real axis. For each $a$, the tip of the slit is located at $-(2-a)/6$, which we can visualize in Fig: \ref{fig: 2a}, \ref{fig: 2b} and \ref{fig: 2c}.\\
		Similarly, proceeding the same way for $F_2(z)=z/(1+z)^2$, where the harmonic shear is given by:
		\begin{align}\label{eqn: 5}
			f_{2,a}(z)=\operatorname{Re}\left\{\frac{a}{3}-\frac{z^2+(a+1)z+\frac{a}{3}+\frac{2}{3}}{(z+1)^3}+\frac{2}{3}\right\}+\dot{\iota}\operatorname{Im}\left\{\frac{z}{(1+z)^2}\right\}.
		\end{align}	 $f(\mathbb{D})$ is the left half plane for a=1 and a slit domain where the tip of the slit is located at $(a+2)/6$ for each $-1\leq a < 1$, as can be seen in Fig: \ref{fig: 2g}, \ref{fig: 2h} and \ref{fig: 2i}.
		Notice that the images of harmonic shear in fig: \ref{fig: 2}, has been changed only horizontally, which is the direction of the shear. We also observe that the shear for $c=2,a=1$ in fig: \ref{fig: 2i} has collapsed the image of unit disk onto the point $1/2$. Therefore, we get the following result:
	\begin{theorem}
		Consider the conformal univalent map $F_{c}$ of the unit disk $\mathbb{D}$ onto  a domain convex in direction of real axis given by (\ref{eqn: 1}), where $c\in[-2,2]$. Let the dilatation function be given by $\omega_{a}(z)=z(z+a)/(1+az)$, where $-1 \leq a \leq 1 $, Then the  horizontal shear of $F_c$ with $c\in(-2,2)$ and dilatation $\omega_{a}$ is given by a two parameter family of harmonic mappings $f_{c,a}(z)=h_{c,a}(z)+\overline{g_{c,a}(z)}$ where $h_{c,a}(z)$ and $g_{c,a}(z)$ are given by (\ref{eqn: 2}) and (\ref{eqn: 3}) such that $f_{c,a}\in S_{H}^0$. For special cases of c=-2 and c=2 with $-1\leq a\leq 1$, the shear is given by (\ref{eqn: 4}) and (\ref{eqn: 5}). Additionally, we have $f(\mathbb{D})$ is a slit domain where the tip of the slit is located at $-(2-a)/6$ for $c=-2$ and at $(a+2)/6$ for $c=2$ such that $		f_{-2,a},f_{2,a}\in S_H^0$.
	\end{theorem}

	\subsection{Minimal surfaces for $a=0$}
	For the special case $a=0$, we have $\omega(z)=z^2$ which is
	a square analytic dilatation. Then we can apply Weierstrass-Enneper formula to lift the harmonic mapping to a minimal graph on $\mathbb{D}$. Thus, we have the following result:
	
	\begin{theorem}
		For $c\in (-2,2)$, define $f_c=h_c+\bar{g_c}:\mathbb{D}\rightarrow \mathbb{C}$ to be the harmonic mapping satisfying $h_{c}(z)-g_{c}(z)=F_{c}(z)$ and $g_{c}'(z)=\omega(z)h_{c}'(z)$ normalized by $h_{c}(0)=g_{c}(0)=g_{c}'(0)=h_{c}'(0)-1=0$, where $F_{c}$ is given by $$F_{c}(z)=\frac{z}{1+cz+z^2},$$ 
		$$\omega(z)=z^2.$$ Then, $f_{c}\in S_{H}^0$ and $f_c(\mathbb{D})$ is convex in the horizontal direction. As $c$ varies from 0 to 2, $f_{c}(\mathbb{D})$ transforms from a infinite vertical strip mapping to a single slit mapping along positive real axis, and as c varies from $-2$ to 0, it transforms from a single slit mapping along negative real axis to a strip mapping. Furthermore, since $w$ is a square analytic function, $f_c$ lifts to a minimal graph $X_{c}$ on $\mathbb{D}$ for each $c\in(-2,2)$ . For some special cases, $c=-2, c=0$ and $c=2, X_{-2}(\mathbb{D})$ is part of Enneper's surface, $X_{0}(\mathbb{D})$ is part of helicoid and $X_{2}(\mathbb{D})$ is again part of Enneper's surface.
	\end{theorem}
		\begin{figure}[H]
		\centering
		\begin{subfigure}[b]{0.49\textwidth}
			\centering
			\includegraphics[width=\textwidth]{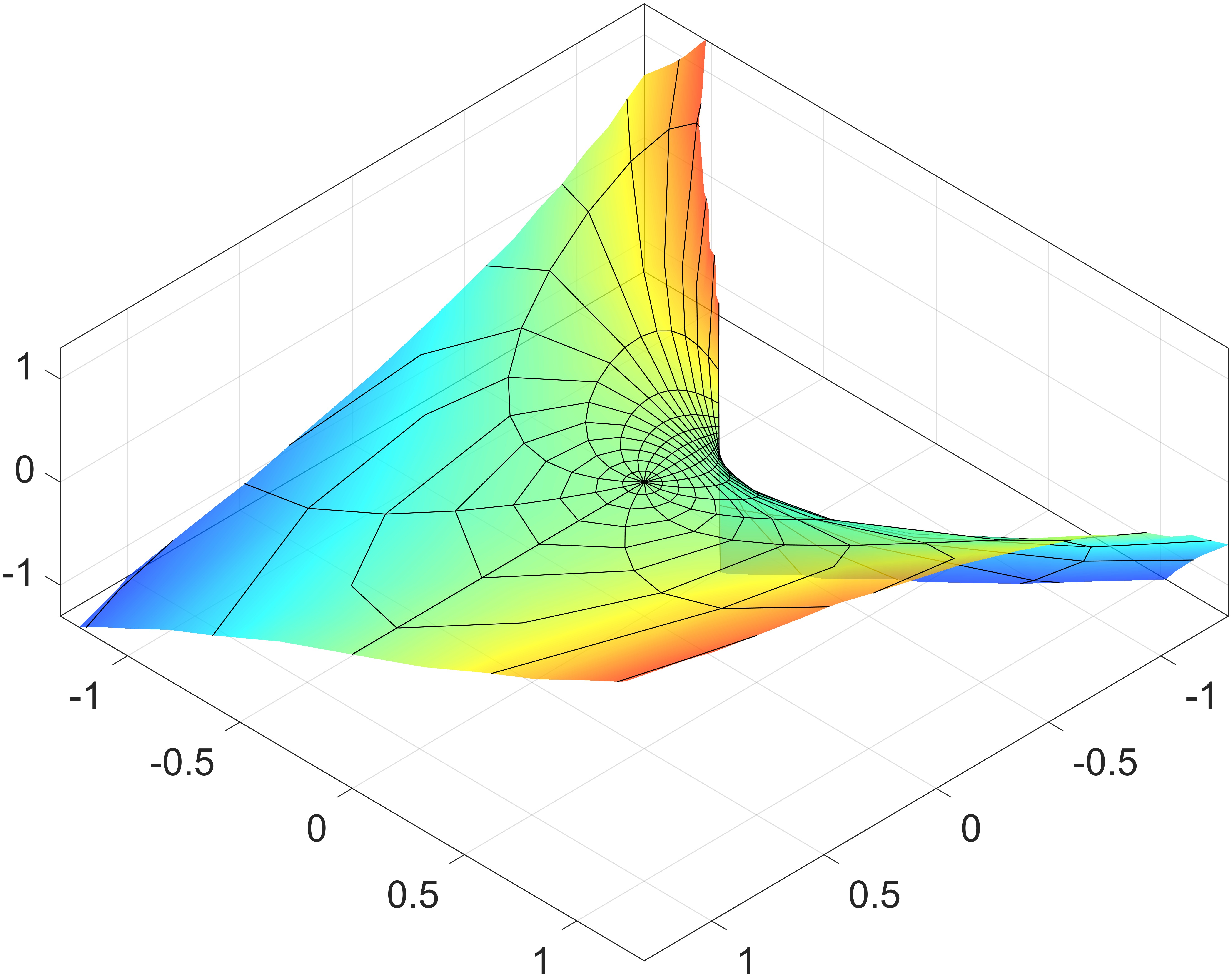}
			\caption{$c=-2$\\ (Identified as Enneper's surface)}
			\label{fig: 3a}
		\end{subfigure}
		\begin{subfigure}[b]{0.49\textwidth}
			\centering
			\includegraphics[width=\textwidth]{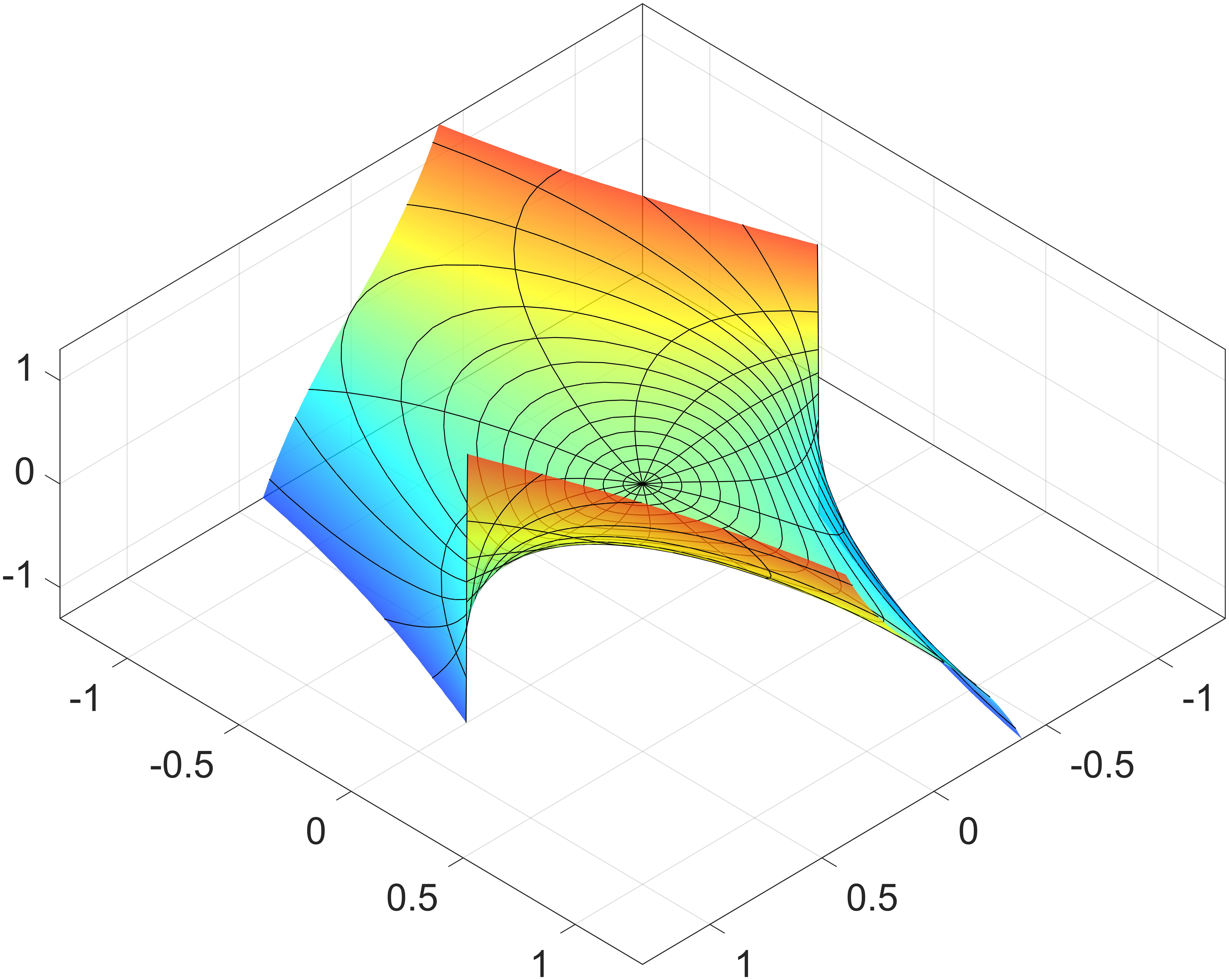}
			\caption{$c=0$\\(Identified as part of Helicoid)}
			\label{fig: 3b}
		\end{subfigure}
		\begin{subfigure}[b]{0.49\textwidth}
			\centering
			\includegraphics[width=\textwidth]{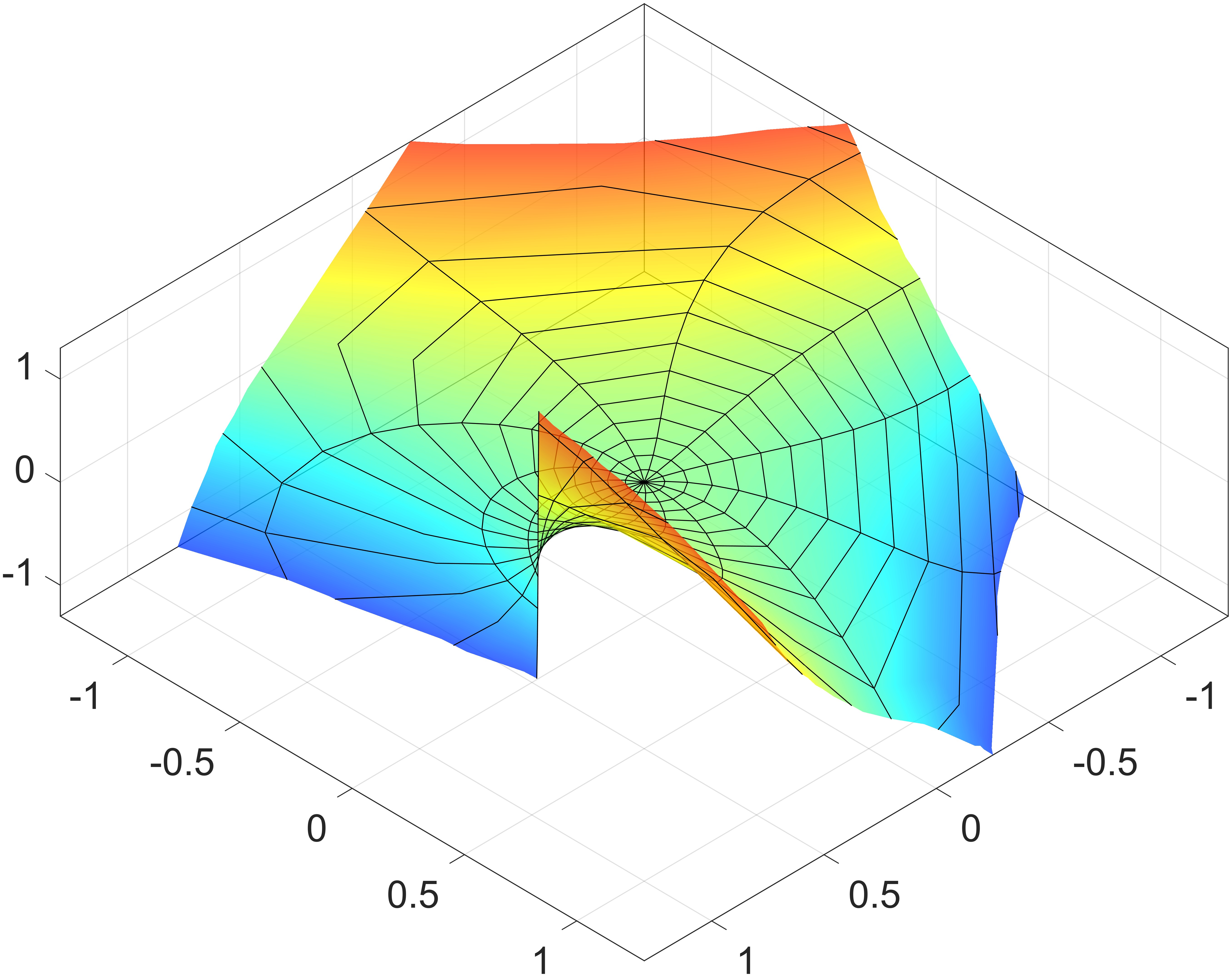}
			\caption{$c=1$\\\;}
			\label{fig: 3c}
		\end{subfigure}
		\begin{subfigure}[b]{0.49\textwidth}
			\centering
			\includegraphics[width=\textwidth]{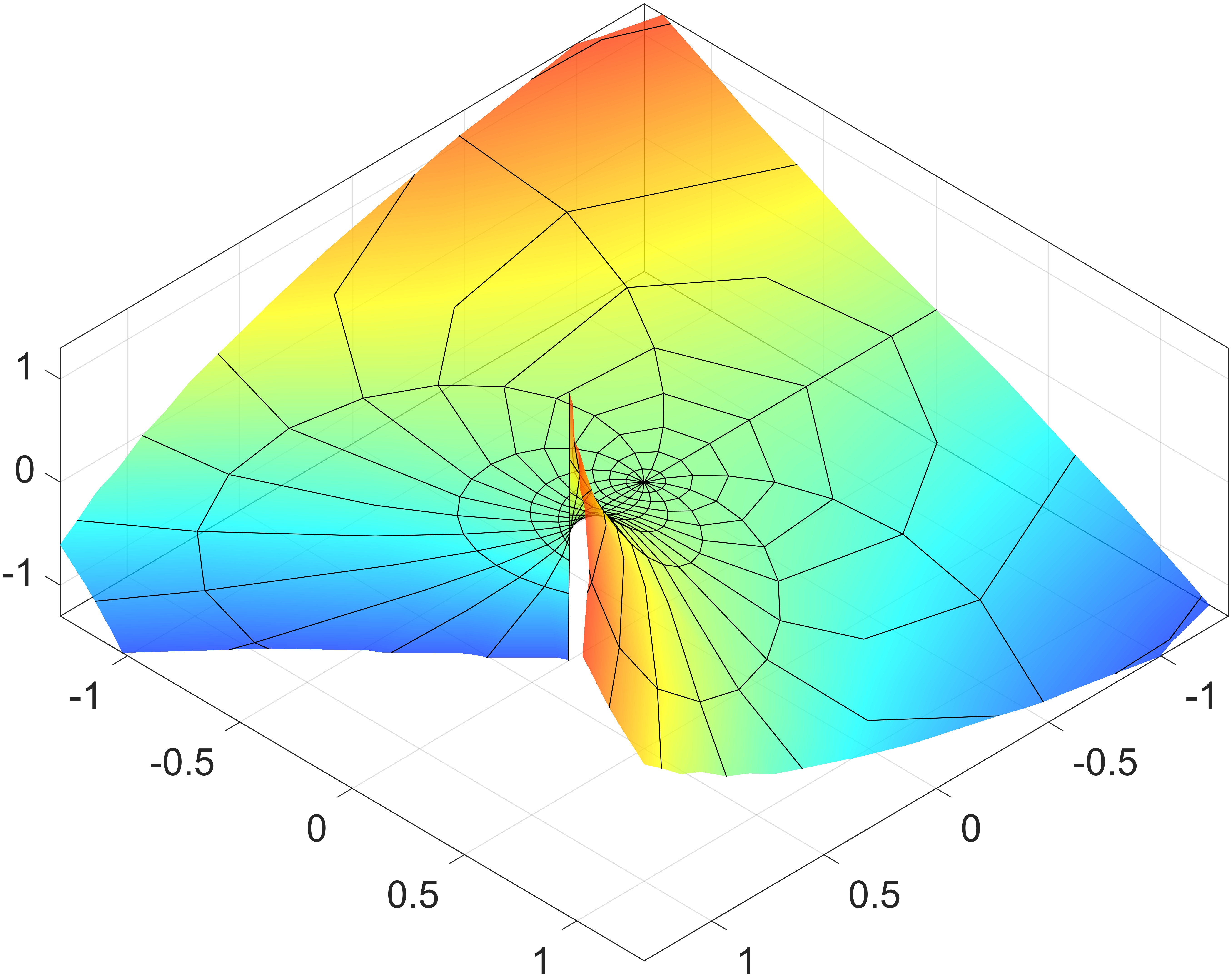}
			\caption{$c=2$\\
				(Identified as Enneper's surface)}
			\label{fig: 3d}
		\end{subfigure}
		\caption{Minimal surfaces over harmonic map $f_c,\, c=0,2,-2,1$ with dilatation $z^2.$}
		\label{fig: 3}
	\end{figure}
	\begin{proof}
		By Theorem \ref{2}, $f_{c}\in S_{H}^0$ and $f_{c}(\mathbb{D})$ is convex in the horizontal direction.Then, the solution of the pair of differential equations $$h_{c}'-g_{c}'=F_{c}'$$
		$$\omega h'-g'=0$$
		defines a family of  harmonic mappings $f_{c}= \operatorname{Re}\{h_{c}+g_{c}\}+\dot{\iota}\operatorname{Im}\{h_{c}-g_{c}\}.$ For $c\in(-2,2)$,
		\begin{align*}
			&u=\operatorname{Re}\{h_{c}+g_{c}\}\\
			&=\operatorname{Re}\left\{\frac{z}{(c-2)(z+1)}+\frac{z}{(c+2)(z-1)}-\frac{8\left(\operatorname{atanh}\left(\frac{c\left(c^2-4\right)}{\sigma_1}\right)-\operatorname{atanh}\left(\frac{\left(c^2-4\right)(c+2z)}{\sigma_1}\right)\right)}{\sigma_1}\right\}\\
			&\,\,\,\,\,\,\,\,\,\,\,-\operatorname{Re}\left\{\frac{z\left(z^2+1\right)}{\left(z^2-1\right)\left(z^2+c z+1\right)}\right\},\\
			&\mbox{where}\,\, \sigma_1=(c-2)^{3 / 2}(c+2)^{3 / 2} \\
			& v=\operatorname{Im}\{h_c-g_c\} =\operatorname{Im}\left\{\frac{z}{z^2+c z+1}\right\} .
		\end{align*}
		Since $\omega(z)=z^2$, by Theorem \ref{2}, $f_{c}$ lifts to a minimal graph corresponding to each $c\in(2,2)$. Applying this theorem yields the following representations of minimal graphs for $c\in(-2,2)$:
		$$X_{c}(z)=(u,v, F(u,v)),$$
		where $x_1$ and $x_2$ are as above, and
		\begin{align*}
			& F(u,v)= 2 \operatorname{Im}\left\{\frac{\frac{2}{c^2-4}+\frac{c z}{c^2-4}}{z^2+c z+1}-\frac{2}{c^2-4}+\frac{2 c \operatorname{atanh}\left(\frac{c\left(c^2-4\right)}{\sigma_1}\right)}{\sigma_1}-\frac{2 c \operatorname{atanh}\left(\frac{\left(c^2-4\right)(c+2 z)}{\sigma_1}\right)}{\sigma_1}\right\} 
		\end{align*}
		where $\sigma_1=(c-2)^{3 / 2}(c+2)^{3 / 2}$.
		We can observe how $f_c(\mathbb{D})$ transforms from an infinite strip to a slit when $c$ goes from 0 to 2 in Fig: \ref{fig: 2e} and \ref{fig: 2h} and vice-versa when $c$ goes from -2 to 0 can be seen in Fig: \ref{fig: 2b} and  Fig: \ref{fig: 2e}.

	 We establish the claim of identifying the surfaces as helicoid and Enneper's surface by considering the following cases:\begin{case}
			$c=-2$, the analytic function becomes $F_{-2}(z)=\dfrac{z}{(z-1)^2}.$
			Applying the shearing technique with dilatation $\omega(z)=z^2$ and then lifitng the resulting harmonic map, we have represetation of the corresponding minimal graph as $X_{-2}(\mathbb{D})=(u,v,F(u,v))$ 
			as:
			\begin{align*}
				X_{-2}(z)=\left(\operatorname{Re}\left\{-\frac{z^2 - z + \frac{2}{3}}{(z -1)^3}-\frac{2}{3}  \right\},\operatorname{Im}\left\{\frac{z}{(z-1)^2}\right\},2\, \operatorname{Im}\left\{\frac{1}{6} - \frac{3z -1}{6(z - 1)^3}\right\}\right).
			\end{align*}
			Applying the change of variables technique to recognize the surface, and substituting $w=(1+z)/(1-z)$ and again expressing as $z=(w-1)/(w+1)$, the surface representation becomes:
			\begin{align*}
				X_{-2}(z)=\left(\operatorname{Re}\left\{\frac{w^3}{12}+\frac{w}{4}-\frac{1}{3}\right\},\operatorname{Im}\left\{\frac{w^2}{4}-\frac{1}{4}\right\},2\, \operatorname{Im}\left\{\frac{w^3}{24} - \frac{w}{8}+\frac{1}{12}\right\}\right).
			\end{align*}
			Interchanging the second and third coordinates, Translation by 1/3, and scaling by a factor of 4, we get 
			\begin{align*}
				X_{-2}(\mathbb{D})=\left(\operatorname{Re}\left\{\frac{w^3}{3}+w\right\}, \operatorname{Im}\left\{\frac{w^3}{3} - w\right\},\operatorname{Im}\left\{w^2\right\} \right).
			\end{align*}
			We see that our original surface $X_{-2}(\mathbb{D})$ is part of the  enneper surface formed by using right half plane as domain instead of the standard unit disk.
		\end{case}
		\begin{case} For c=0,
			the function $F_{c}(z)$ becomes
			 $F(z)=z/(1+z^2)$  with dilatation  $\omega(z)=z^2.$ 
			Solving for $h$ and $g$, we get
			 $h(z)=(z+\mathrm{atan}\left(z\right)+z^2 \,\mathrm{atan}\left(z\right))/(2\,{\left(z^2 +1\right)})$ and 
				$g(z)=(\mathrm{atan}\left(z\right)-z+z^2 \,\mathrm{atan}\left(z\right))/(2\,{\left(z^2 +1\right)})$, such that $f=h+\overline{g}$
			is the harmonic map required to lift to a minimal graph. So, we have
			\begin{align*}
				X_{0}(z)=(x_1,x_2,x_3)
				&=\left(\operatorname{Re}\{\mathrm{atan}\left(z\right)\},\operatorname{Im}\{\frac{z}{z^2 +1}\},2\operatorname{Im}\left\{\frac{1}{2}-\frac{1}{2\,{\left(z^2 +1\right)}}\right\}\right)\\
				&=\left(\operatorname{Re}\left\{\frac{\dot{\iota}}{2}\,\log\left(\frac{\dot{\iota}+z}{\dot{\iota}-z}\right)\right\},\operatorname{Im}\left\{\frac{z}{z^2 +1}\right\},2\operatorname{Im}\left\{\frac{1}{2}-\frac{1}{2\,{\left(z^2 +1\right)}}\right\}\right),
			\end{align*}
			which are hardly recognizable as representations of part of
			the helicoid.
			However, applying the change
			of variables $w =(\dot{\iota}+z)/(\dot{\iota}-z)$ and again writing 
			$z=\dot{\iota}(w-1)/(w+1)$,
			\begin{align*}
				X_{0}(z)&=\left(\operatorname{Re}\left\{\frac{\dot{\iota}}{2}\,\log\left(w\right)\right\},\operatorname{Im}\left\{\frac{\dot{\iota}}{4}\left(w-\frac{1}{w}\right)\right\},2\operatorname{Im}\left\{-\frac{1}{8}\left(w+\frac{1}{w}-2\right)\right\}\right).
			\end{align*}
			Scaling by 4 gives,
			\begin{align*}
				X_{0}(z)
				&=\left(\operatorname{Im}\{2\log\left(w\right)\},\operatorname{Re}\left\{w-\frac{1}{w}\right\},\operatorname{Im}\left\{-\left(w+\frac{1}{w}\right)\right\}\right).
			\end{align*}
			Since, negative sign over $x_3$, interchanging of the second and third coordinates, and then first and third coordinates won't alter the geometry and gives us the formula of helicoid. So we see $X_{0}(\mathbb{D})$ is the same surface as $Y_0({z\in \mathbb{C}: \operatorname{Re}{z}>0}))$, i.e., $X_0(\mathbb{D})$ is part of the helicoid.
		\end{case}

		\begin{case}
			For c=2, the conformal map takes the form $F_2(z)=\dfrac{z}{(1+z)^2}.$
			After similar computations we get,
			\begin{align*}
				X_{2}(z)=\left(\operatorname{Re}\left\{\frac{2}{3} - \frac{z^2 + z + \frac{2}{3}}{(z + 1)^3}\right\},\operatorname{Im}\left\{\frac{z}{(1+z)^2}\right\},2\, \operatorname{Im}\left\{\frac{1}{6} - \frac{\frac{z}{2} + \frac{1}{6}}{(z + 1)^3}\right\}\right).
			\end{align*}
			Substituting $w=(1-z)/(1+z)$ and expressing in terms of $z$ as  $z=(1-w)/(1+w)$,then translating by $-1/3$, scaling by a factor of $-1/4$ and finally  adjusting the negative over $x_3$  the surface representation becomes:
			 \begin{align*}
			 	X_{2}(z)=\left(\operatorname{Re}\left\{\frac{w^3}{3}+w\right\},\operatorname{Im}\left\{w^2\right\}, \operatorname{Im}\left\{\frac{w^3}{3} -w\right\}\right).
			 \end{align*}
			 We see interchanging the coordinates suitably gives the parameterization of Enneper's surface.\\
		\end{case}
		We have plotted the minimal surfaces formed for $c=-2, 0, 1$ and $2$ in Fig: \ref{fig: 3a}, \ref{fig: 3b}, \ref{fig: 3c} and \ref{fig: 3d}. Their projection onto $x$-$y$ plane can be observed and matched with images of the corresponding harmonic maps in Fig: \ref{fig: 2b}, \ref{fig: 2e} and \ref{fig: 2h}.
	\end{proof}

	Further, this idea of lifting can be generalized by changing the dilatation to $\omega(z)=z^{2n}$, $n\in \mathbb{N}$. In 2014, Ponnusamy et al. have constructed a two-parameter family of harmonic mappings $f_{n,c}$ with dilatation $z^n$ by shearing the four slit conformal mapping  defined by $\phi(z)$ in \cite{ponnusamy2014harmonic}. We will further lift the resulting shears to a family of minimal surfaces for even $n$ and $A=0$.
	\begin{align*}
	\phi(z)=A\,\log\left(\frac{1+z}{1-z}\right)+B\, \frac{z}{1+cz+z^2}.
	\end{align*}
	Taking $A=0, B=1$ serves our purpose. By writing $c=-2\,\cos\gamma$ where $\gamma \in (0,\pi)$. Applying shearing technique gives, $$h(z)=-\frac{\dot{\iota}}{2\,\sin\gamma}(e^{-\dot{\iota}\gamma}\textit{I}_2-e^{\dot{\iota}\gamma}\textit{I}_3),$$
	\begin{align*}
		\text{where} \,\,\textit{I}_2=\int_{0}^{z}\frac{d\zeta}{(\zeta-e^{-\dot{\iota}\gamma})^2(1-\zeta^{n})}\,\,\, \mbox{and}\,\,\,\textit{I}_3=\int_{0}^{z}\frac{d\zeta}{(\zeta-e^{\dot{\iota}\gamma})^2(1-\zeta^{n})}.
	\end{align*}
	Using $\dfrac{1}{1-z^n}=-\frac{1}{n}\sum \limits_{k=0}^{n-1}\dfrac{z_k}{z-z_k}$, where $z_k=e^{\frac{2\pi\dot{\iota}k}{n}}$ for $k=0,1,2,...n-1$ and assuming that $\eta=e^{\dot{\iota}\gamma} \neq z_k$, consider
	\begin{align*}
		\textit{I}_{\eta}&=\int_{0}^{z}\frac{d\zeta}{(\zeta-\eta)^2(1-\zeta^{n})}\\
		&=-\frac{1}{n}\sum \limits_{k=0}^{n-1}\left\{\frac{1}{\eta-z_k}\left(\frac{1}{\eta-z}-\frac{1}{\eta}\right)-\frac{1}{(\eta-z_k)^2}\left[\log\left(\frac{\eta-z}{\eta}\right)-\log\left(\frac{z_k-z}{z_k}\right)\right]\right\}.
	\end{align*}
	 \begin{figure}
		\centering
		\begin{subfigure}[b]{0.49\textwidth}
			\centering
			\includegraphics[width=\textwidth]{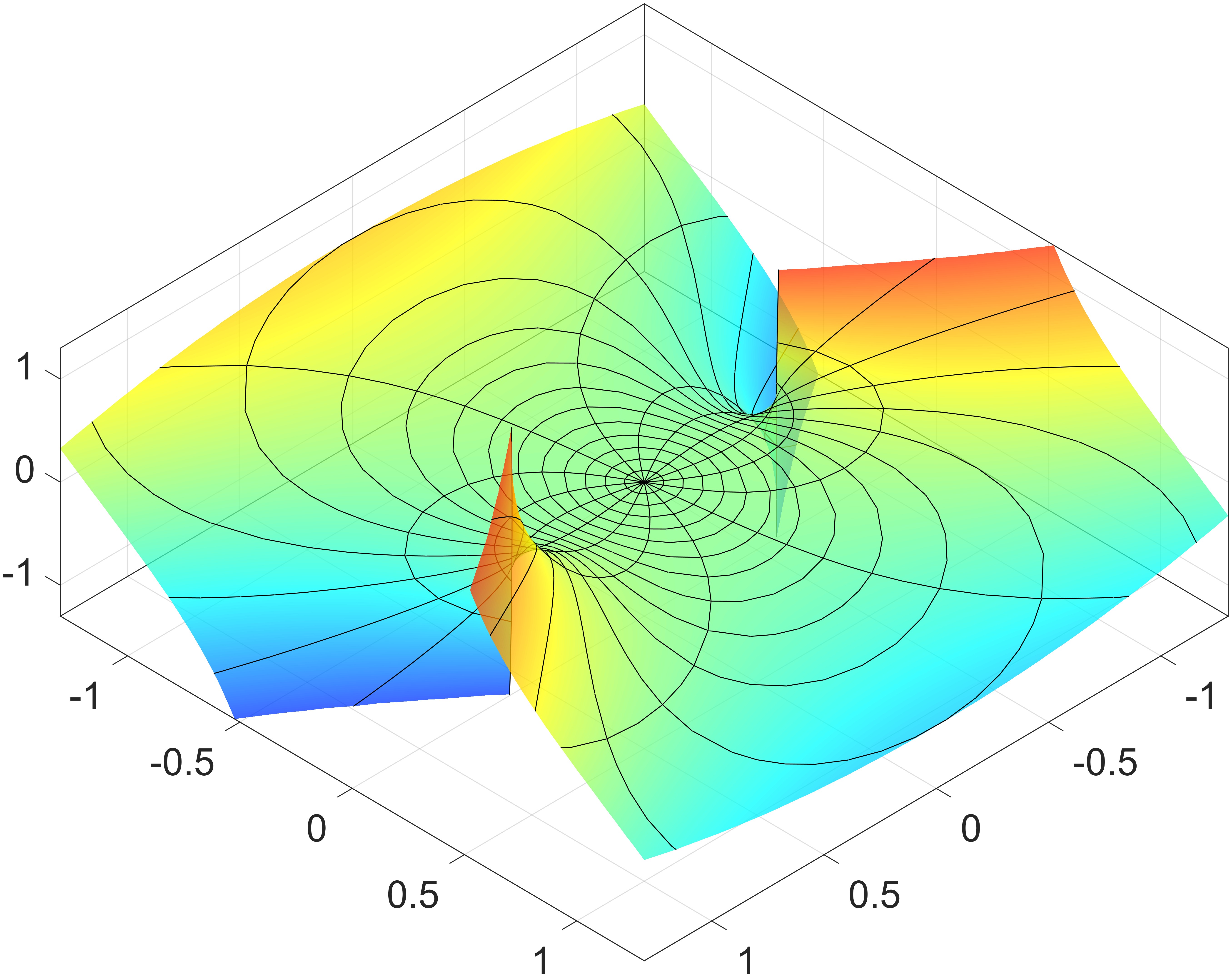}
			\caption{c=0, n=4}
			\label{fig: 4a}
		\end{subfigure}
		\begin{subfigure}[b]{0.49\textwidth}
			\centering
			\includegraphics[width=\textwidth]{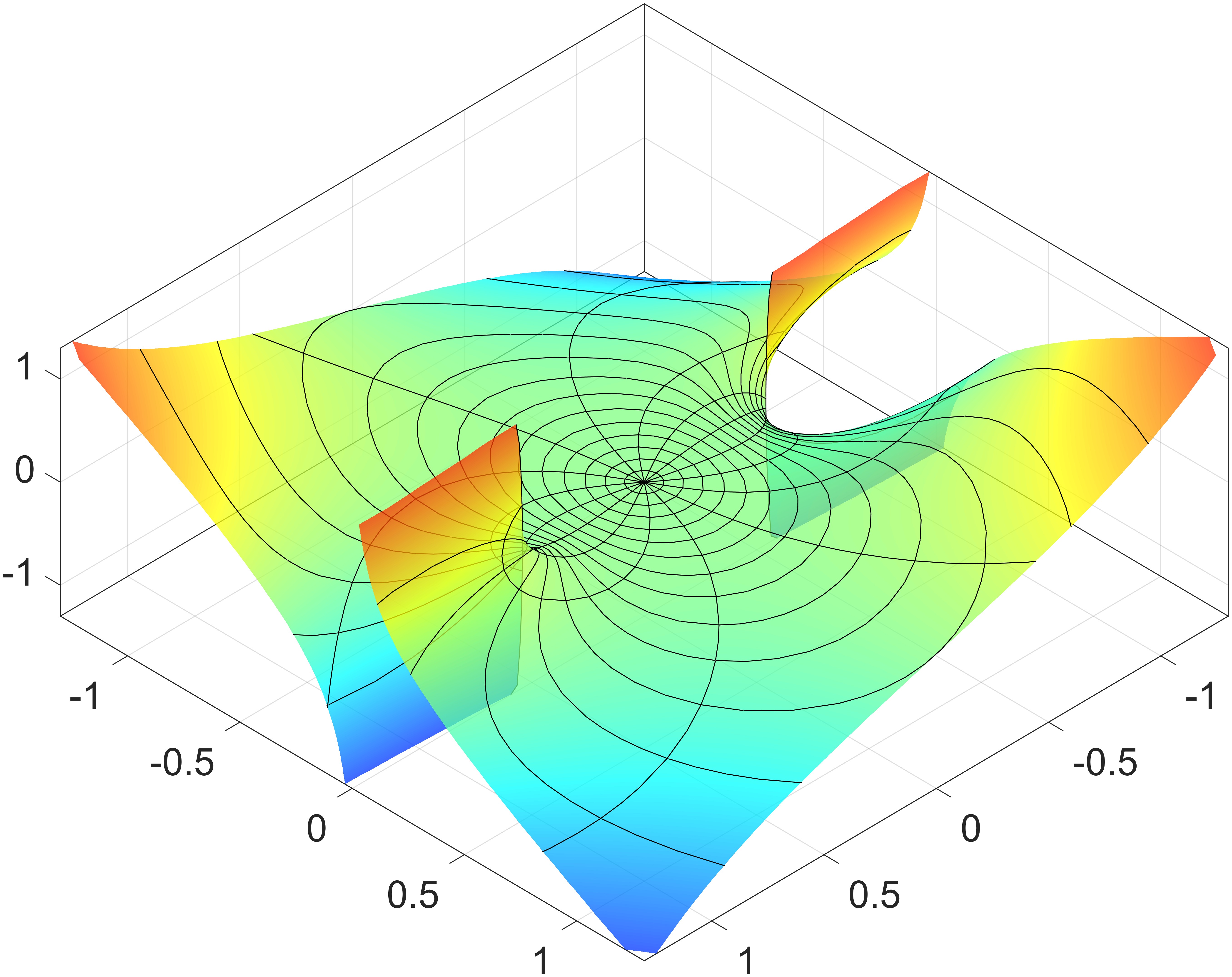}
			\caption{$c=0, n=6$}
			\label{fig: 4b}
		\end{subfigure}
		\begin{subfigure}[b]{0.49\textwidth}
			\centering
			\includegraphics[width=\textwidth]{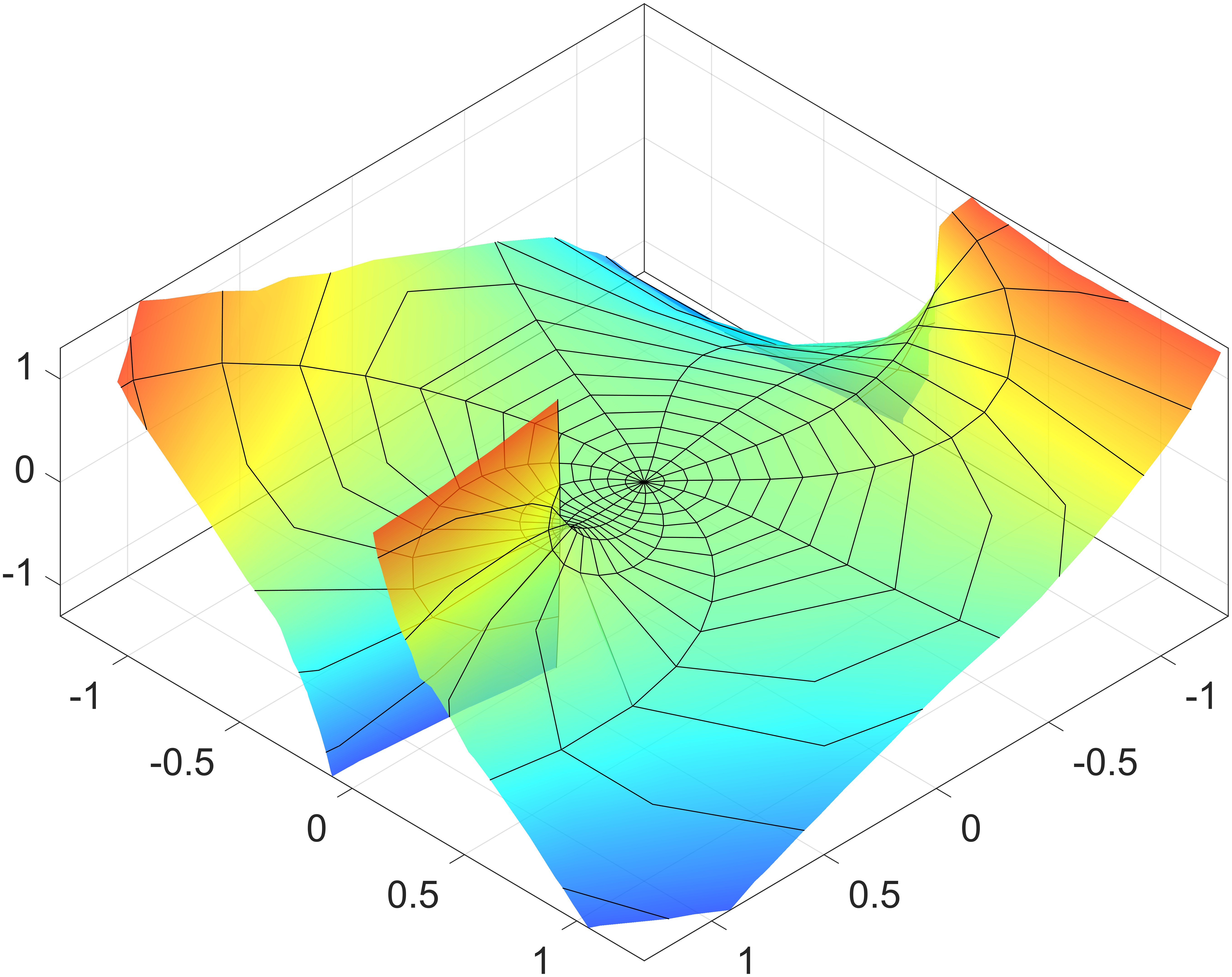}
			\caption{$c=1, n=4$}
			\label{fig: 4c}
		\end{subfigure}
		\begin{subfigure}[b]{0.49\textwidth}
			\centering
			\includegraphics[width=\textwidth]{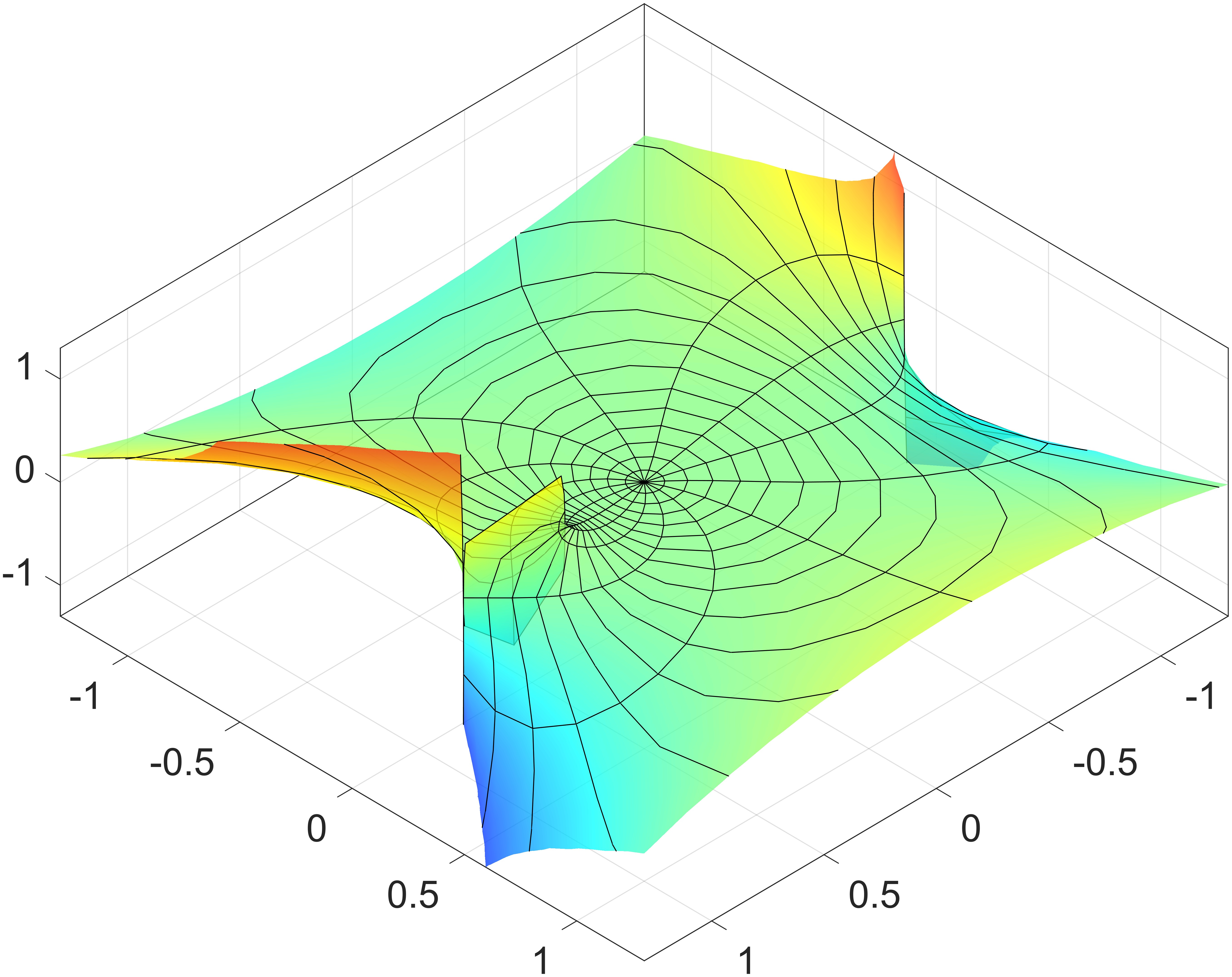}
			\caption{$c=1, n=6$}
			\label{fig: 4d}
		\end{subfigure}
		\begin{subfigure}[b]{0.49\textwidth}
			\centering
			\includegraphics[width=\textwidth]{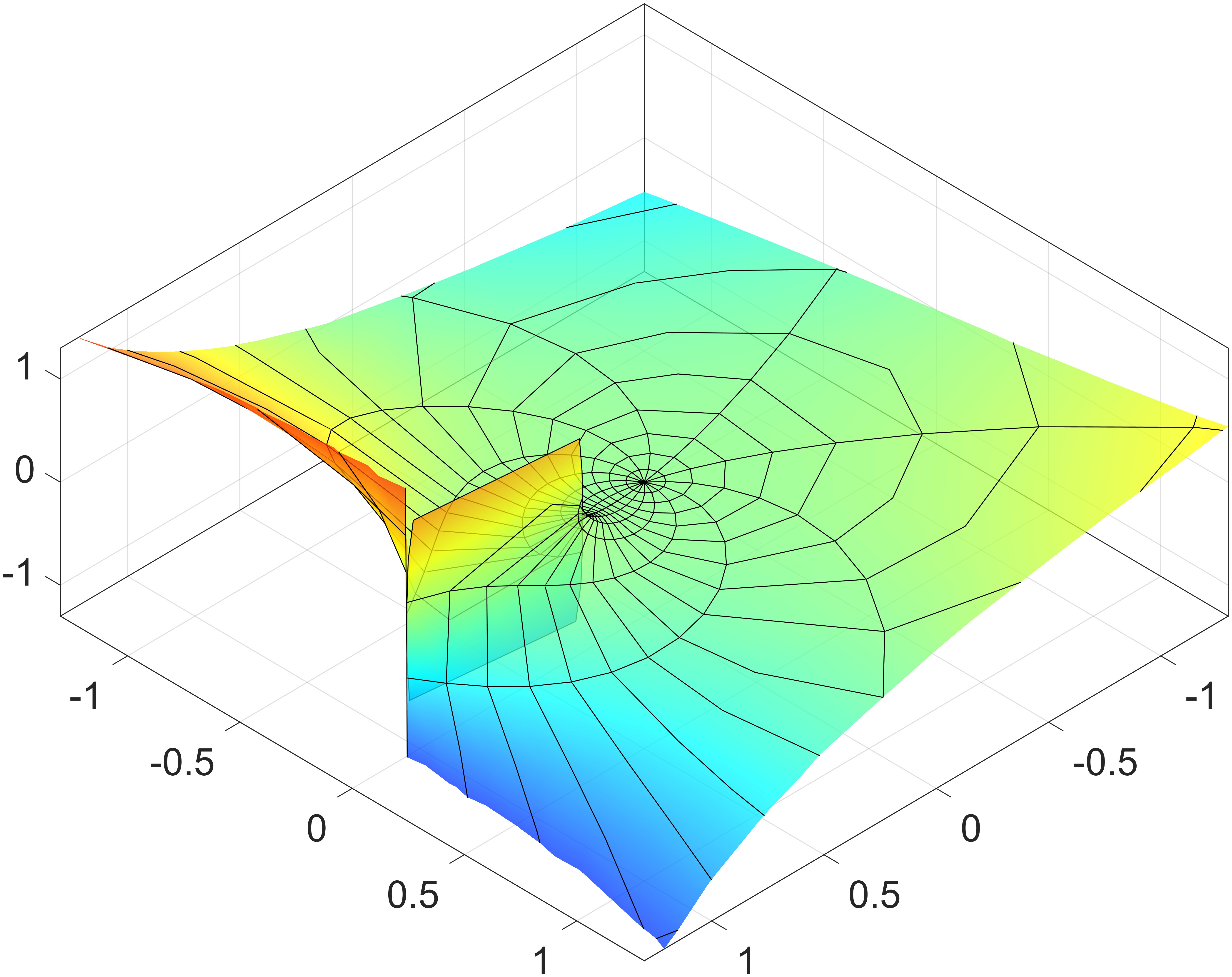}
			\caption{$c=2, n=4$}
			\label{fig: 4e}
		\end{subfigure}
		\begin{subfigure}[b]{0.49\textwidth}
			\centering
			\includegraphics[width=\textwidth]{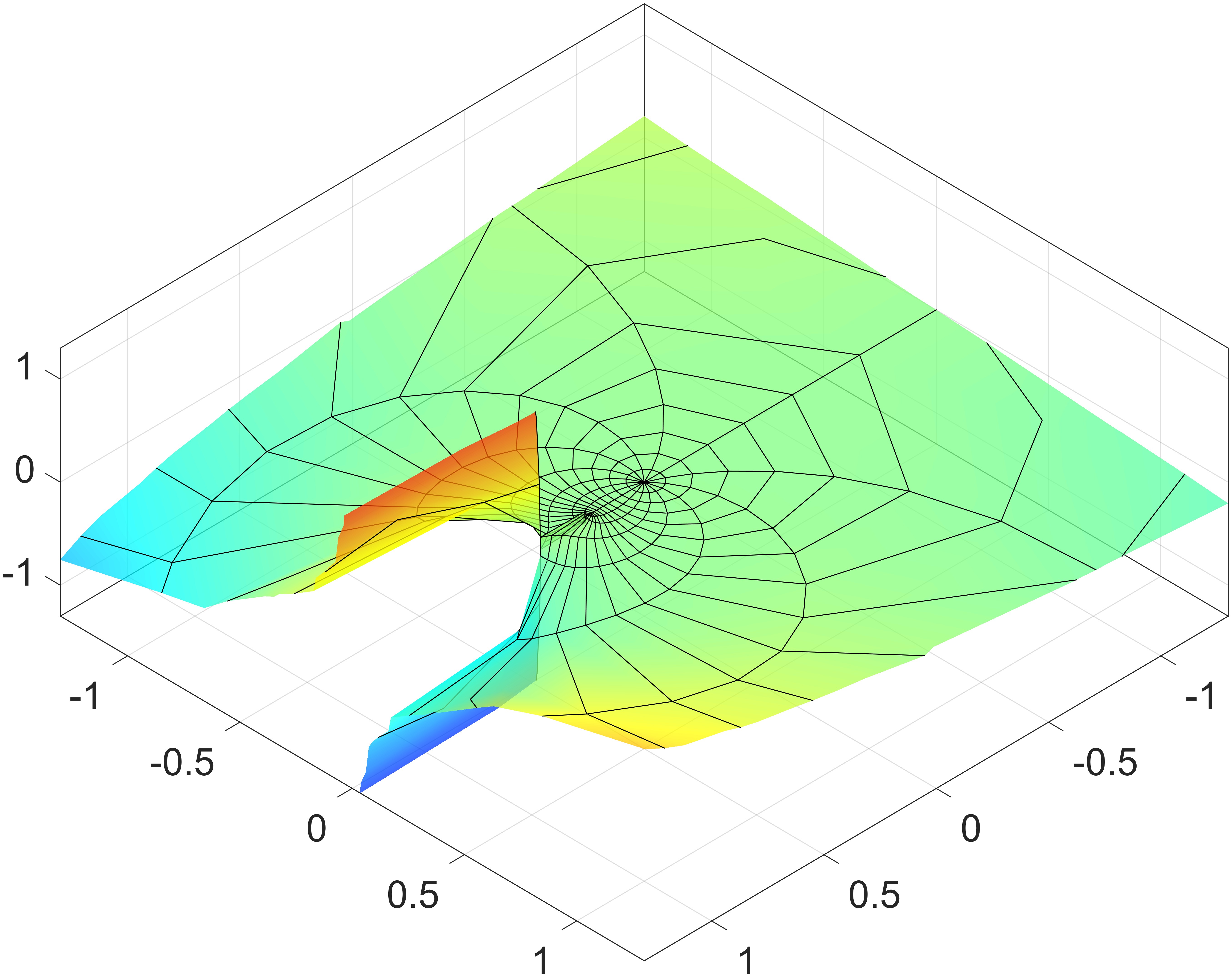}
			\caption{$c=2, n=6$}
			\label{fig: 4f}
		\end{subfigure}
		\caption{Minimal surfaces over harmonic shears of $F(z)=z/(1+cz+z^2)$ with dilatation $\omega(z)=z^n$ for even n.}
		\label{fig: 4}
	\end{figure}
	Let $\mathbf{N}=\{0,1,\dots,n-1\}$ be an index set. Suppose that $a\in \mathbf{N}$, then we define $\mathbf{N}_a=\mathbf{N}\setminus\{a\}$.\\
	In the case of $\gamma=2\pi m/n$, where $m=0,1,2,\dots,n-1$, 
	\begin{align*}
		\textit{I}_{3,m}&=\int_{0}^{z}\frac{d\zeta}{(\zeta-z_m)^2(1-\zeta^{n})}\\
		&=\frac{1}{n}\sum \limits_{k \in \mathbf{N}_m}\int_{0}^{z}\frac{z_k\,d\zeta}{(\zeta-z_m)^2(\zeta-z_k)}+\frac{1}{n}\int_{0}^{z}\frac{z_m\,d\zeta}{(\zeta-z_m)^3}.
	\end{align*}
	The above sum can be computed as above and the last integral is,
	\begin{align*}
		\frac{1}{n}\int_{0}^{z}\frac{d\zeta}{(\zeta-z_m)^3}=\frac{1}{2n}\left[\frac{1}{(z-z_m)^2}-\frac{1}{z_m^2}\right].
	\end{align*}
     We can readily find $g$ from $h-\phi$ for $A=0$. Thus, we get the harmonic map $h=f+\overline{g}$ that we can lift to a minimal surface for square dilatation by taking $\omega(z)=z^{2n}$. Since, $x_3=2\operatorname{Im}\left\{\int_{0}^{z}\sqrt{\omega(\zeta)}\,h'(\zeta)\, d\zeta\right\}$, we compute
    $$\int_{0}^{z}\sqrt{\omega(\zeta)}\,h'(\zeta)\, d\zeta=  -\frac{\dot{\iota}}{2\,\sin\gamma}\left(e^{-\dot{\iota}\gamma}\int_{0}^{z}\frac{\zeta^n }{(\zeta-e^{-\dot{\iota}\gamma})^2(1-\zeta^{2n})}d\zeta-e^{\dot{\iota}\gamma} \int_{0}^{z}\frac{\zeta^n}{(\zeta-e^{\dot{\iota}\gamma})^2(1-\zeta^{2n})} d\zeta \right).$$
     Using, $\dfrac{z^n}{1-z^{2n}}=\dfrac{z^n}{(1-z^n)(1+z^n)}=\dfrac{1}{2}\left[\dfrac{1}{1-z^n}-\dfrac{1}{1+z^n}\right],$
   the above expression becomes:
     \begin{align*}
     &= \frac{-\dot{\iota}}{2\,\sin\gamma}\left(e^{-\dot{\iota}\gamma}\int_{0}^{z}\frac{1 }{2(\zeta-e^{-\dot{\iota}\gamma})^2(1-\zeta^n)}-e^{-\dot{\iota}\gamma}\int_{0}^{z}\frac{1 }{2(\zeta-e^{-\dot{\iota}\gamma})^2(1+\zeta^n)}\right) d\zeta\\
     &\,\,\,\,\,\,-\frac{-\dot{\iota}}{2\,\sin\gamma}\left(e^{\dot{\iota}\gamma} \int_{0}^{z}\frac{1}{2(\zeta-e^{\dot{\iota}\gamma})^2(1-\zeta^n)} -e^{\dot{\iota}\gamma} \int_{0}^{z}\frac{1}{2(\zeta-e^{\dot{\iota}\gamma})^2(1+\zeta^n)}\right) d\zeta.
     \end{align*}
     We observe that the first and third integral in above expression are $\textit{I}_2$ and $\textit{I}_3.$ For the second (say $\textit{I}_4$) and fourth integral ($\textit{I}_5$), again by partial fractions we can express:\\
     $$\frac{1}{1+z^n}=-\frac{1}{n}\sum \limits_{l=0}^{n-1}\frac{z_l}{z-z_l},$$ 
     where $z_l$'s are $n^{th}$ roots of equation $z^n+1=0$ such that
      $z_l=e^{\frac{2\pi\dot{\iota}l}{n}+\frac{\pi \dot{\iota}}{n}}$, where $l=0,1,2,\dots,n-1$, then
      $$\int_{0}^{z}\sqrt{\omega(\zeta)}\,h'(\zeta)\, d\zeta= \frac{-\dot{\iota}}{2\,\sin\gamma}\left(\frac{e^{-\dot{\iota}\gamma}}{2}\textit{I}_2-\frac{e^{-\dot{\iota}\gamma}}{2}\textit{I}_4-\frac{e^{\dot{\iota}\gamma}}{2}\textit{I}_3+\frac{e^{\dot{\iota}\gamma}}{2}\textit{I}_5\right),$$
     where,
      \begin{align*}
      	\textit{I}_4=\int_{0}^{z}\frac{1 }{(\zeta-e^{-\dot{\iota}\gamma})^2(1+\zeta^n)} d\zeta \,\,\,\mbox{and}\,\,\,\textit{I}_5=\int_{0}^{z}\frac{1 }{(\zeta-e^{\dot{\iota}\gamma})^2(1+\zeta^n)} d\zeta.
    \end{align*}
  Assuming that $\rho=e^{\dot{\iota}\gamma} \neq z_l$, i.e, $\gamma \neq (2 l+1)\frac{\pi}{n}$, and computing in a similar manner just like $\textit{I}_{\eta}$, we get
      \begin{align*}
      	\textit{I}_{\rho}&=\int_{0}^{z}\frac{d\zeta}{(\zeta-\rho)^2(1-\zeta^{n})}\\
      	&=-\frac{1}{n}\sum \limits_{k=0}^{n-1}\left\{\frac{1}{\rho-z_l}\left(\frac{1}{\rho-z}-\frac{1}{\rho}\right)-\frac{1}{(\rho-z_l)^2}\left[\log\left(\frac{\rho-z}{\rho}\right)-\log\left(\frac{z_l-z}{z_l}\right)\right]\right\}.
      \end{align*} 
     In case of $\gamma = (2 s+1)\frac{\pi}{n}$ for $s=0,1,2,\dots,n-1$ and computing in similar manner just like  $\textit{I}_{3,m}$ above, we get  $\textit{I}_{3,p}=\textit{I}_5$.
     \begin{align*}
     	\textit{I}_{3,s}&=\int_{0}^{z}\frac{d\zeta}{(\zeta-z_s)^2(1-\zeta^{n})}=-\frac{1}{n}\sum \limits_{l \in \mathbf{N}_s}\int_{0}^{z}\frac{z_l\,d\zeta}{(\zeta-z_s)^2(\zeta-z_l)}+\frac{1}{n}\int_{0}^{z}\frac{z_s\,d\zeta}{(\zeta-z_s)^3}.
     \end{align*}
     This last integral can again be computed as done for $I_{3,m}$ above. In Fig: \ref{fig: 4}, we have plotted the minimal surfaces obtained after calculating $x_3$ coordinate, for different values of $c$ and even $n$.
     
	 \begin{table}[ht!]
	 	\caption{The integrals $I_2$ and $I_3$ for the analytic part $h=\frac{1}{2 \sin \gamma} i\left(e^{-i \gamma} I_2-e^{i \gamma} I_3\right)$ of the harmonic shear $f$ with a dilatation $\omega(z)=z^n$.}
	 	\label{Tab: 1}
	 	\centering
	 	\begin{tabular}{lll}
	 		\hline$\gamma$ & $I_2$ & $I_3$ \\
	 		\hline is not $2 \pi m / n$ & $I_{\bar{\eta}}$ & $I_{\eta}$ \\
	 		is $2 \pi m / n$ & $I_{3, n-m}$ & $I_{3, \mathrm{~m}}$ \\
	 		\hline
	 	\end{tabular}
	 \end{table}
	 \begin{table}[ht!]
	 	\caption{The integrals $I_4$ and $I_5$ for the $x_3=\operatorname{Im}\left\{\frac{-\dot{\iota}}{2\,\sin\gamma}\left(e^{-\dot{\iota}\gamma}2\textit{I}_2-e^{-\dot{\iota}\gamma}\textit{I}_4-e^{\dot{\iota}\gamma}\textit{I}_3+e^{\dot{\iota}\gamma}\textit{I}_5\right)\right\}$ coordinate for the minimal surfaces.}
	 	\label{Tab: 2}
	 	\centering
	 	\begin{tabular}{lll}
	 		\hline$\gamma$ & $I_4$ & $I_5$ \\
	 		\hline is not $(2s+1) \pi / n\,\,\,\,\,$ & $I_{\bar{\rho}}$ & $I_{\rho}$ \\
	 		is $(2s+1) \pi / n$ & $I_{3, n-s}$ & $I_{3,s}$ \\
	 		\hline
	 	\end{tabular}
	 \end{table}
	 \section{Harmonic shears of inner region of an epicycloid}
\begin{figure}[H]
	 	\centering
	 	\begin{subfigure}[b]{\textwidth}
	 		\centering
	 		\includegraphics[width=\textwidth]{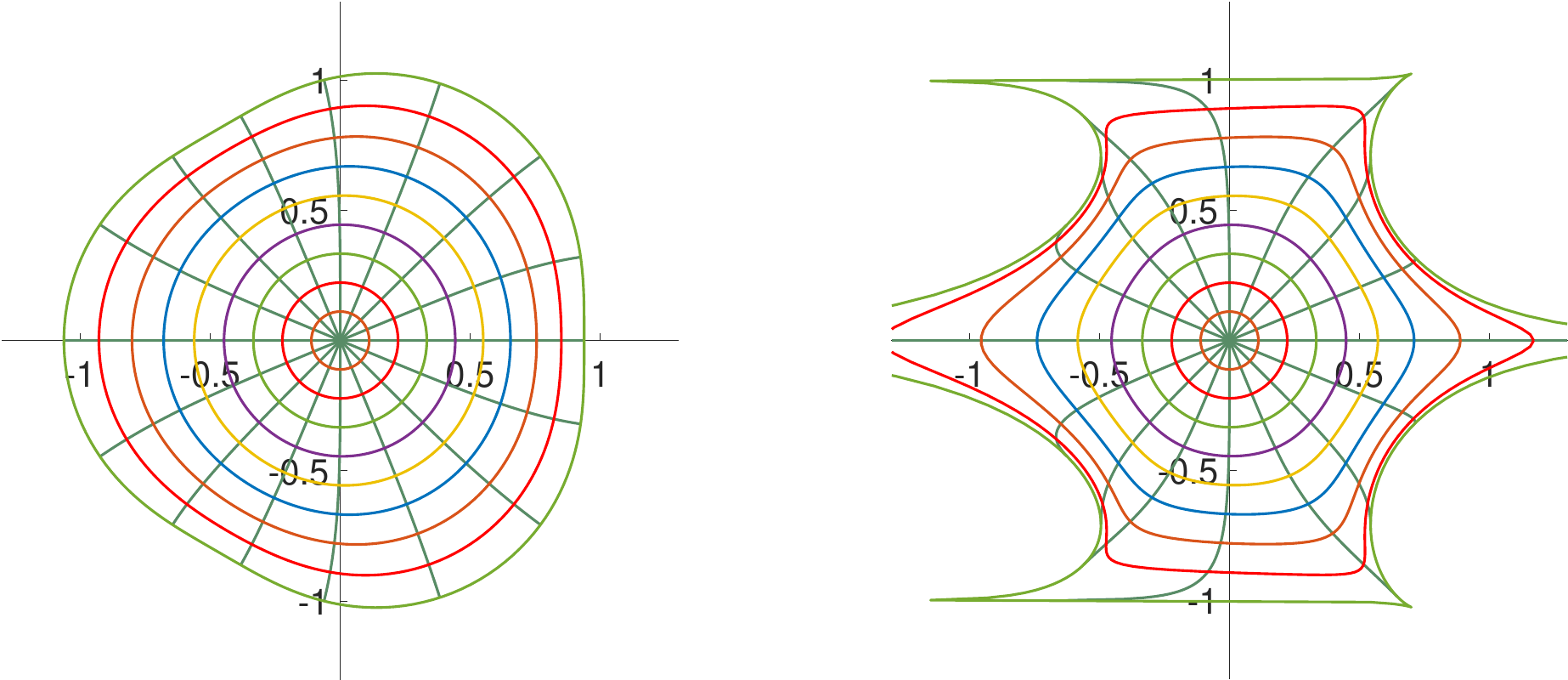}
	 		\caption{$n=2$}
	 		\label{fig: 5a}
	 	\end{subfigure}
	 	\begin{subfigure}[b]{\textwidth}
	 		\centering
	 		\includegraphics[width=\textwidth]{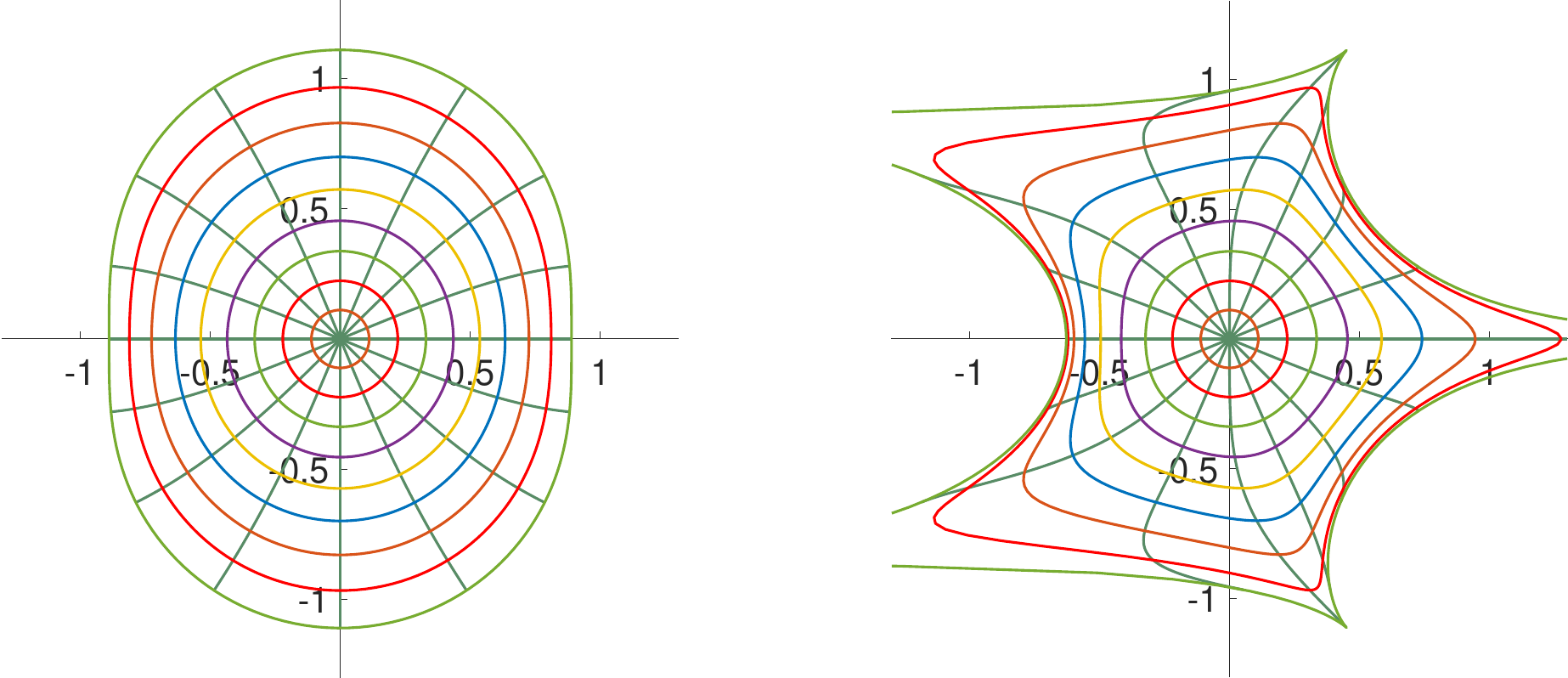}
	 		\caption{$n=3$}
	 		\label{fig: 5b}
	 	\end{subfigure}
	 	\begin{subfigure}[b]{\textwidth}
	 		\centering
	 		\includegraphics[width=\textwidth]{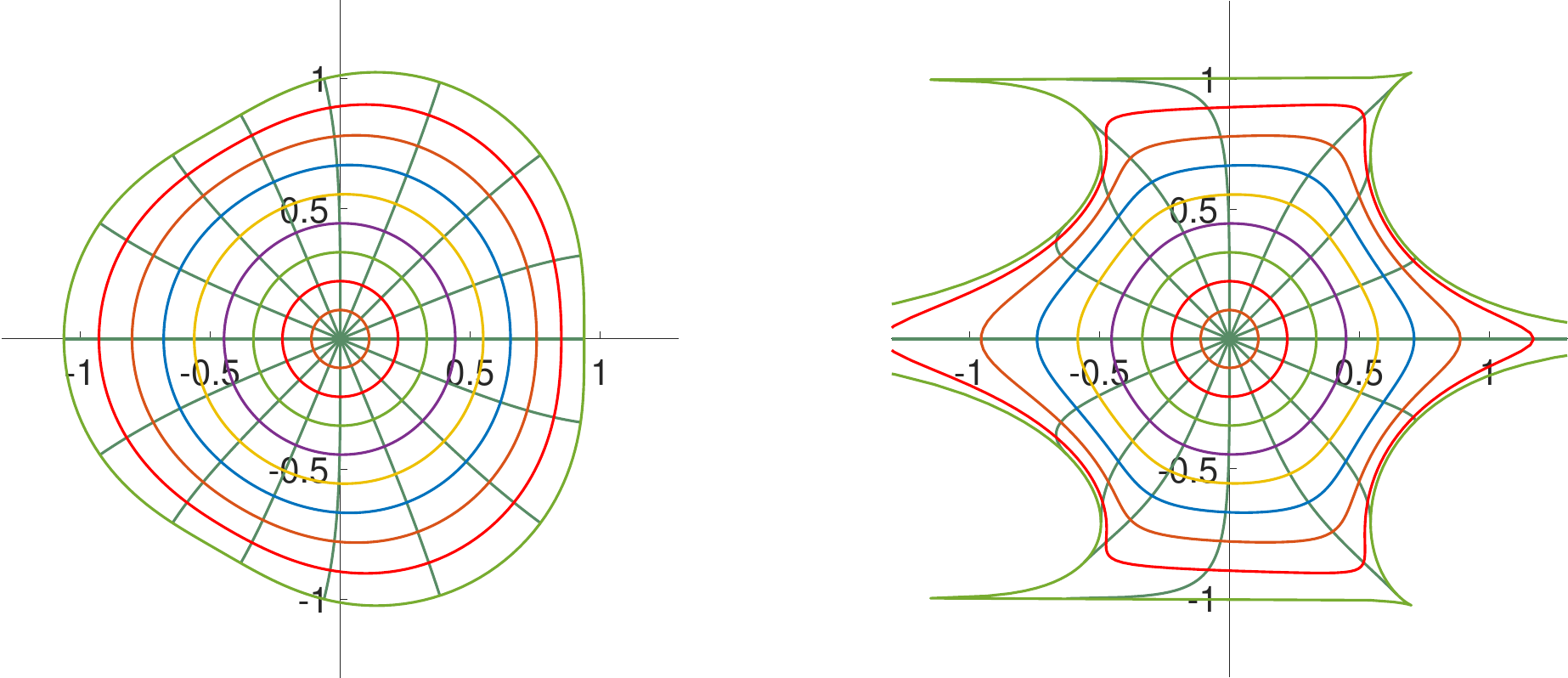}
	 		\caption{$n=4$}
	 		\label{fig: 5c}
	 	\end{subfigure}
	 	\caption{Conformal mapping $F_n$ of the unit disk onto inner region of an \textit{epicycloid} with $n-1$ cusps, alongside its harmonic shear $f_n$ on right with $n+2$ concave boundary arcs for the dilatation $z^n,\, n=2,3,4$.}
	 	\label{fig: 5}
	 \end{figure}

\begin{figure}
	 	\centering
	 	\begin{subfigure}[b]{0.5\textwidth}
	 		\centering
	 		\includegraphics[width=\textwidth]{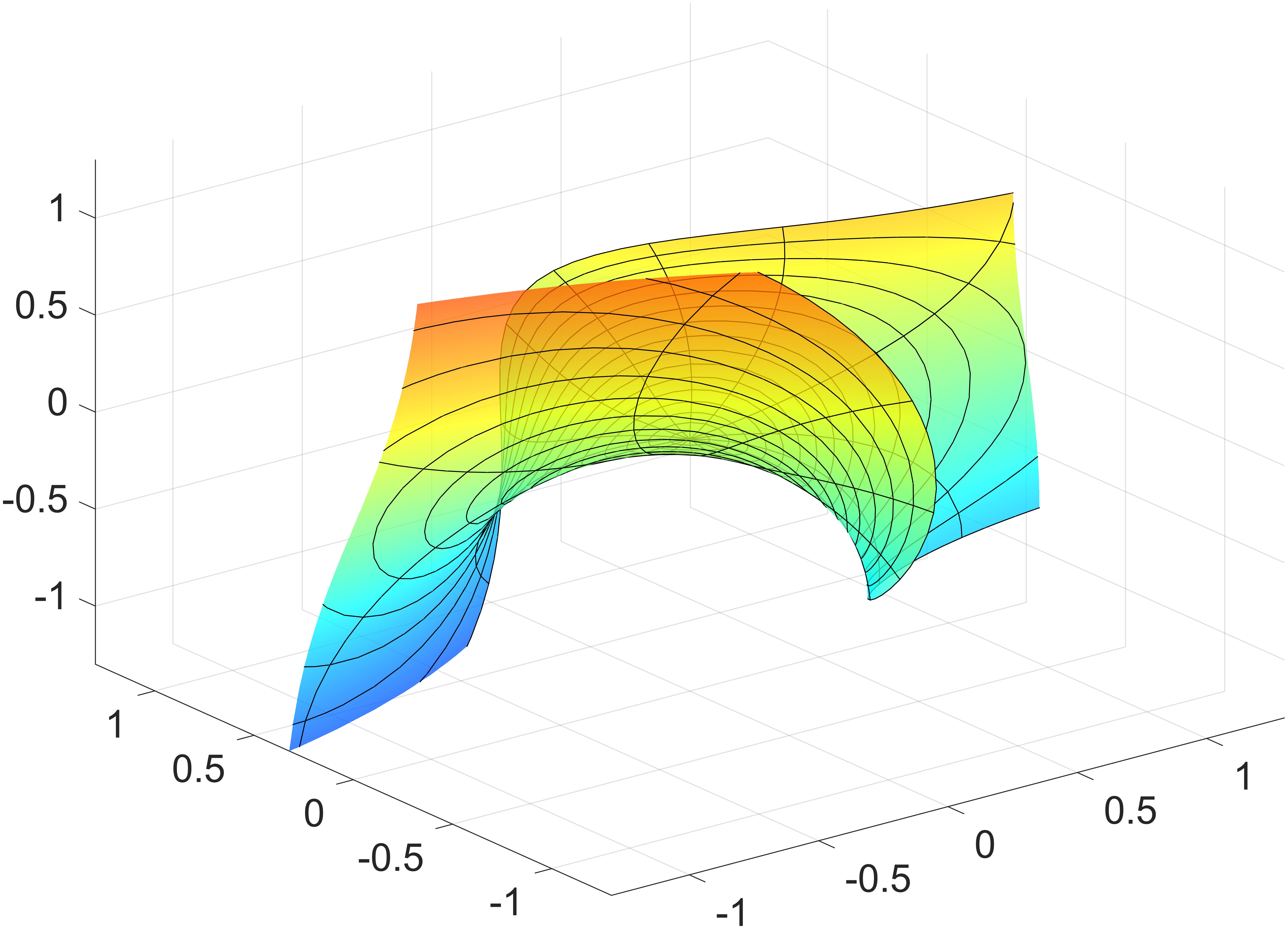}
	 		\caption{$n=2$.}
	 		\label{fig: 6a}
	 	\end{subfigure}
	 	\begin{subfigure}[b]{0.5\textwidth}
	 		\centering
	 		\includegraphics[width=\textwidth]{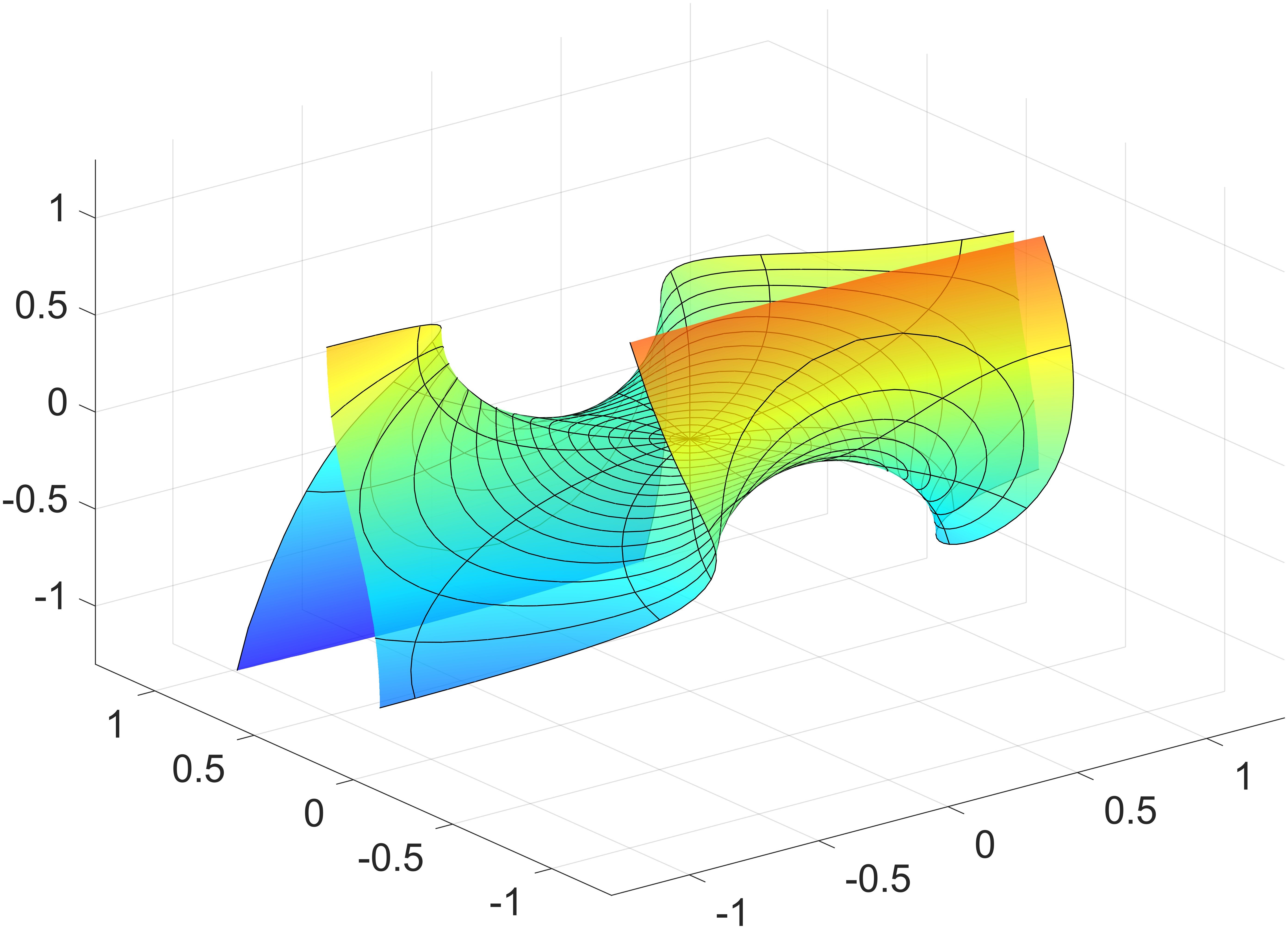}
	 		\caption{$n=3$.}
	 		\label{fig: 6b}
	 	\end{subfigure}
	 	\begin{subfigure}[b]{0.5\textwidth}
	 		\centering
	 		\includegraphics[width=\textwidth]{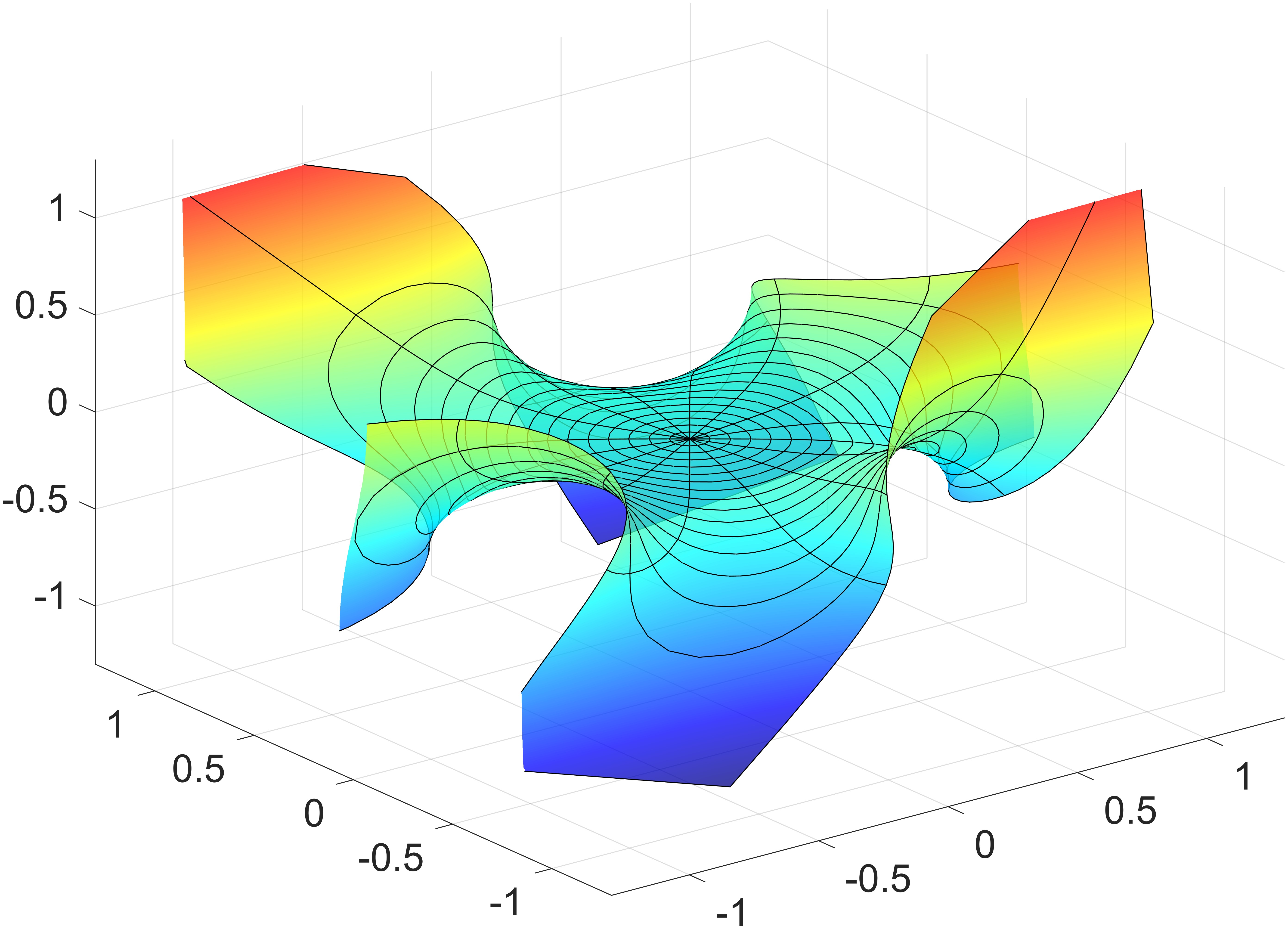}
	 		\caption{$n=4$.}
	 		\label{fig: 6c}
	 	\end{subfigure}
	 	\caption{Minimal surface over harmonic shear $f_n$ with dilatation $\omega=z^{2n}$ for $n=2,3,4.$}
	 	\label{fig: 6}
	 \end{figure}

	 In this section, we will apply shear construction method on a new class of conformal mappings denoted by $F_n(z)$ with dilatation $\omega(z)=z^n$ and will further lift to a family of minimal surfaces using Theorem \ref{3}. We have created images of these minimal surfaces for $n=2, 3, 4$ in Fig: \ref{fig: 6}. We have also plotted images of conformal mappings $F_n(z)$ alongside images of corresponding harmonic mappings $f_n(z)$ for $n=2, 3, 4$ in Fig: \ref{fig: 5}.\\
	 let 
	 \begin{align}\label{eqn: 6}
	 	F_{n}(z)=h_{n}(z)-g_{n}(z)=z-\frac{1}{n^2}\, z^n,
	 \end{align}
	 be the class of conformal univalent mappings of the unit disk $\mathbb{D}$ onto  a domain convex in horizontal direction and $n$ is any natural number greater than equal to $2$. We observed that images of the open unit disk under conformal maps form the interior region of an \textit{epicycloid} characterized by $n-1$ cusps.  Let
	 \begin{align*}
	 	\omega(z)=\frac{g'(z)}{h'(z)}=z^n,
	 \end{align*}
	 be the dilatation. Differentiating $(1)$ gives us, 
	 \begin{align*}
	 	h_{n}'(z)-g_{n}'(z)=1-\frac{1}{n}\cdot z^{n-1}.
	 \end{align*}
	 Considering the pair of these differential equations above, we get
	 $$h_{n}'(z)=\frac{1-\frac{1}{n}\cdot z^{n-1}}{1-z^n}.$$
	 Integrating from $0$ to $z$
	 and normalizing such that $h(0)=0$ and similarly, solving for $g(z)$ gives: 
	 \begin{align}
	 	& h_{n}(z)= z \hspace{0.08cm} _2F_1\left(1, \frac{1}{n} ; \frac{1}{n}+1 ; z^n\right)+\frac{\log \left(z^n-1\right)}{n^2}-\frac{\pi \dot{\iota}}{n^2},\label{eqn: 7} \\
	 	& g_{n}(z)= z \hspace{0.08cm} _2F_1\left(1, \frac{1}{n} ; \frac{n+1}{n} ; z^n\right)-z+\frac{\log \left(z^n-1\right)+z^n}{n^2} -\frac{\pi \dot{\iota}}{n^2} \label{eqn: 8},
	 \end{align}
	 where ${ }_2 F_1(a, b ; c ; z)$ represents a hypergeometric function, which is a power series such that,
	 $$
	 { }_2 F_1(a, b ; c ; z)=\sum_{n=0}^{\infty} \frac{(a)_n(b)_n}{(c)_n} \frac{z^n}{n !},
	 $$
	 where $a, b$, and $c \in \mathbb{C}$ and
	 $
	 (x)_n=x(x+1) \cdots(x+n-1)
	 $
	 is the Pochhammer symbol. So, the desired map corresponding to analytic map $F_n$ is $f_n=h_n+\overline{g_n}$. By Theorem \ref{2}, $f_n \in S_H^0$ and $f_n$ is convex in horizontal direction.
	 \par We see that image of equally spaced radial segments and concentric circles under sense-preserving harmonic univalent map $f_n$ is a radial plot with n+2 concave boundary arcs. Notice that, although the image of $F_n$ is bounded, the image of $f_n$ is unbounded.The harmonic map $f_n$ can also be expressed as $f_n(z)=x_1(z)+\dot{\iota}x_2(z)$ where
	 \begin{align*}
	 	x_{1}(z)&=\operatorname{Re}\left\{h(z)+g(z)\right\} \\
	 	&=\operatorname{Re}\left\{2z \hspace{0.08cm} _2F_1\left(1, \frac{1}{n} ; \frac{n+1}{n} ; z^n\right)-z+\frac{2\log \left(z^n-1\right)+z^n}{n^2} -\frac{2\pi \dot{\iota}}{n^2}\right\} ,\\
	 	x_{2}(z)&=\operatorname{Im}\{h(z)-g(z)\}\\
	 	&= z-\frac{z^n}{n^2}.
	 \end{align*}
	 
	 For $n=2m$,\,$f_n$ lifts to a minimal graph where 
	 \begin{align*}
	 	x_3&=2\,\operatorname{Im}\left\{ \int_{0}^{z}\sqrt{\omega(\zeta)}h_{n}'(\zeta)\,d\zeta\right\}\\
	 	&=2\,\operatorname{Im}\left\{ \int_{0}^{z}\zeta^m\left(\frac{1-\frac{1}{2m}\cdot \zeta^{2m-1}}{1-\zeta^{2m}}\right)
	 	\,d\zeta\right\}\\
	 	&=2\,\operatorname{Im}\left\{ \int_{0}^{z}\left(\frac{\zeta^m}{1-\zeta^{2m}}-\frac{1}{2m}\cdot\frac{ \zeta^{3m}}{\zeta(1-\zeta^{2m})}\right)
	 	\,d\zeta\right\}\\
	 	&=2\,\operatorname{Im}\left\{ \int_{0}^{z}\left(\frac{1}{2}\left(\frac{1}{1-\zeta^{m}}-\frac{1}{1+\zeta^{m}}\right)
	 	-\frac{1}{2m}\cdot\left(\frac{\zeta^{m-1}}{2(\zeta^{m}+1)}-\frac{\zeta^{m-1}}{2(\zeta^{m}-1)}-\zeta^{m-1}\right)\right)
	 	\,d\zeta\right\}\\
	 	&=2\,\operatorname{Im}\left\{z \hspace{0.08cm} _2F_1\left(1, \frac{1}{m} ; \frac{1}{m}+1 ; z^m\right)-z \hspace{0.08cm} _2F_1\left(1, \frac{1}{m} ; \frac{1}{m}+1 ; -z^m\right)-\frac{1}{2m}\left(\frac{\mathrm{atanh}\left(z^m \right)-z^m }{m}\right)\right\}.\,\,(9)
	 \end{align*}
	 Thus, we obtain the following result:
	 \begin{theorem}
	 	Consider the conformal univalent map of the unit disk $\mathbb{D}$ onto  a domain convex in direction of real axis given by $(\ref{eqn: 6})$. Let the dilatation function be given by $\omega(z)=z^n$ ,then the  horizontal shear of $F_n$ with dilatation $\omega$ is given by $f_n(z)=h_n(z)+\overline{g_n(z)}$, where $h_{n}(z)$ and $g_{n}(z)$ are given by (\ref{eqn: 7}) and (\ref{eqn: 8}) such that $f_{n}\in S_{H}^0$. In particular, when $n$ is an even positive integer, $f_n(\mathbb{D})$ lifts to a minimal graph $(x_1,x_2,x_3)$ where $x_1=\operatorname{Re}\left\{h(z)+g(z)\right\}$, $x_2=\operatorname{Im}\{h(z)-g(z)\}$ and $x_3$ is given by $(9)$.
	 	
	 \end{theorem}
	 The images of the unit disk under $F_{n}$ and $f_{n}$ for $n=2,3$ and $4$ are shown in Fig: \ref{fig: 5} as plots of the images of equally spaced radial segments and concentric circles. The images of the resulting minimal surfaces obtained by lifting $f_n$ in $\mathbb{R}^3$ are depicted in Fig: \ref{fig: 6}.
	 
	 \bibliographystyle{ieeetr}

	 \end{document}